\documentclass[aos,MSNbibl,nameyear,seceqn,dvips]{arximspdf}
\usepackage{graphicx}
\doi{10.1214/12-AOS1079} 
\volume{41}
\issue{1}
\pubyear{2013}
\firstpage{143}
\lastpage{176}

\makeatletter

\newcommand{\iint}{\int\!\!\int}

\newtheorem{theorem}{Theorem}[section]
\newtheorem{lemma}{Lemma}[section]
\newtheorem{corollary}{Corollary}[section]

\newproclaim{example}{Example}[section]
\newproclaim{remark}{Remark}[section]

\newcommand{\cal}{\mathcal}

\newcommand{\ve}{\varepsilon}
\newcommand{\Real}{\mathbb R}
\newcommand{\dd}{\mathrm{d}}
\newcommand{\tr}{\operatorname{tr}}

\newcommand{\bvarphi}{\bolds{\varphi}}

\newcommand{\bM}{{M}}
\newcommand{\bB}{{B}}
\newcommand{\bC}{{C}}
\newcommand{\bD}{{D}}
\newcommand{\bff}{{f}}

\makeatother

\begin{document}
\begin{frontmatter}

\title{Optimal design for linear models with
correlated~observations\thanksref{T1}}
\runtitle{Optimal design for correlated observations}

\thankstext{T1}{Supported in part by the Collaborative Research
Center ``Statistical modeling of nonlinear dynamic processes'' (SFB
823, Teilprojekt C2) of the German Research Foundation (DFG).}

\begin{aug}
\author[A]{\fnms{Holger}~\snm{Dette}\corref{}\ead[label=e1]{holger.dette@rub.de}},~%
\author[B]{\fnms{Andrey}~\snm{Pepelyshev}\thanksref{t2}\ead[label=e2]{pepelyshev@stochastik.rwth-aachen.de}}
\and
\author[C]{\fnms{Anatoly}~\snm{Zhigljavsky}\ead[label=e3]{ZhigljavskyAA@cf.ac.uk}}
\runauthor{H. Dette, A. Pepelyshev and A. Zhigljavsky}
\affiliation{Ruhr-Universit\"at Bochum, RWTH Aachen and Cardiff University}
\address[A]{H. Dette\\
Fakult\"at f\"ur Mathematik\\
Ruhr-Universit\"at Bochum\\
Bochum, 44780\\
Germany\\
\printead{e1}}
\address[B]{A. Pepelyshev\\
Institute of Statistics\\ RWTH Aachen University\\
Aachen, 52056\\
Germany\\
\printead{e2}}
\address[C]{A. Zhigljavsky\\
School of Mathematics\\
Cardiff University\\
Cardiff, CF24 4AG\\
United Kingdom\\
\printead{e3}} 
\end{aug}

\thankstext{t2}{Supported in part by Russian Foundation of Basic
Research, Project 12-01-007470.}

\received{\smonth{4} \syear{2012}}
\revised{\smonth{9} \syear{2012}}

%
\begin{abstract}
In the common linear regression model the problem of determining
optimal designs for least squares estimation is considered in the
case where the observations are correlated. A necessary condition for
the optimality of a given design is provided, which extends the
classical equivalence theory for optimal designs in models with
uncorrelated errors to the case of dependent data. If the regression
functions are eigenfunctions of an integral operator defined by the
covariance kernel, it is shown that the corresponding measure defines a
universally optimal design.
For several models universally optimal designs can be identified
explicitly. In particular, it is proved that the uniform distribution
is universally optimal for a class of trigonometric regression models
with a broad class of covariance kernels and that the arcsine
distribution is universally optimal for the polynomial regression model
with correlation structure defined by the logarithmic potential. To
the best knowledge of the authors these findings provide the first
explicit results on optimal designs for regression models with
correlated observations, which are not restricted to the location scale
model.
\end{abstract}

%
\begin{keyword}[class=AMS]
\kwd[Primary ]{62K05}
\kwd{60G10}
\kwd[; secondary ]{31A10}
\kwd{45C05}
\end{keyword}
\begin{keyword}
\kwd{Optimal design}
\kwd{correlated observations}
\kwd{integral operator}
\kwd{eigenfunctions}
\kwd{arcsine distribution}
\kwd{logarithmic potential}
\end{keyword}
\end{frontmatter}

\section{Introduction}
\label{intro}

Consider the common linear regression model
%
%
\begin{equation}
\label{lin-model} y(x)=\theta_1f_1(x)+\cdots+
\theta_mf_m(x)+\ve(x),
\end{equation}
where $f_1(x),\ldots, f_m(x)$ are linearly independent, continuous functions,
$\ve(x)$ denotes a random error process or field, $\theta_1,\ldots,
\theta_m$ are unknown parameters and $x$ is the explanatory
variable,\vadjust{\goodbreak}
which varies in a compact design space
$\mathcal{X}\subset\mathbb{R}^d$. We assume that $N$ observations, say
$y_1,\ldots, y_N$, can be taken at experimental conditions
$x_1,\ldots,x_N$ to estimate the parameters in the linear regression
model (\ref{lin-model}). If an appropriate estimate, say $\hat\theta$,
of the parameter $\theta=(\theta_1,\ldots,\theta_m)^T$ has been chosen,
the quality of the statistical analysis can be further improved by
choosing an appropriate design for the experiment. In particular, an
optimal design minimizes a functional of the variance--covariance matrix
of the estimate $\hat\theta$, where the functional should reflect
certain aspects of the goal of the experiment. In contrast to the case
of uncorrelated errors, where numerous results and a rather complete
theory are available [see, e.g., the monograph of
\citet{pukelsheim2006}],
the construction of optimal designs for dependent observations is
intrinsically more difficult.
On the other hand, this problem is of particular practical interest as
in most applications there exists correlation between different observations.
Typical examples include models, where the explanatory variable $x$
represents the time and all observations correspond to one subject.
In such situations optimal experimental designs are very difficult to
find even in simple cases. Some exact optimal design problems were
considered in \citet{bolnat1982}, N{\"a}ther [(\citeyear{naether1985a}),
Chapter 4], N{\"a}ther
(\citeyear{naether1985b}), \citet{pazmue2001} and \citet{muepaz2003},
who derived optimal designs for the location scale model
%
%
\begin{equation}
\label{const} y(x) = \theta+ \varepsilon(x).
\end{equation}
Exact optimal designs for specific linear models have been investigated
in \citet{detkunpep2008,kiessteh2008,harstu2010}.
Because explicit solutions of optimal design problems for correlated
observations are rarely available, several authors have proposed to determine
optimal designs based on asymptotic
arguments [see, e.g., Sacks and Ylvisaker (\citeyear{sackylv1966,sackylv1968}), \citet
{bickherz1979}, \citet{naether1985a}, \citet{zhidetpep2010}].
Roughly speaking, there exist three approaches to embed the optimal
design problem for regression models with correlated observations in an
asymptotic
optimal design problem. The first one is due to
Sacks and Ylvisaker (\citeyear{sackylv1966,sackylv1968}), who assumed that
the covariance structure of the error process $\ve(x)$ is fixed and
that the number of design points tends to infinity.
Alternatively, \citet{bickherz1979}
and \citet{bickherzsch1981} considered a
different model, where the correlation function depends on the sample size.
Recently, \citet{zhidetpep2010} extended the Bickel--Herzberg approach
and allowed the variance (in addition to the correlation function) to vary
as the number of observations changes. As a result,
the corresponding optimality criteria contain a kernel with singularity
at zero.
The focus in all these papers is again mainly on the location scale
model~(\ref{const}).

The difficulties in the development of the optimal
design theory for correlated observations can be explained by
a different structure of the covariance of the least squares estimator
in model (\ref{lin-model}), which is of the form $M^{-1} B M^{-1}$
for certain matrices $M$ and $B$ depending on the design.
As a consequence, the corresponding design problems are in general not
convex [except for
the location scale model
(\ref{const}) where $M=1$].

The present paper is devoted to the problem of determining optimal designs
for more general models with correlated observations than the
simple location scale model (\ref{const}).
In Section~\ref{prelim} we present some preliminary discussion and
introduce the necessary notation.
In Section~\ref{equiv} we investigate general conditions for design optimality.
One of the main results of the paper is Theorem~\ref{thuniversal-optimality},
where we derive necessary and sufficient conditions for the universal
optimality of designs.
By relating the optimal design problems to eigenvalue problems
for integral operators we identify a broad class of multi-parameter
regression models
where the universally optimal designs can be determined explicitly. It
is also shown that in this case the least squares estimate with the
corresponding optimal design has the same covariance matrix as the
weighted least squares estimates with its optimal design. In other
words, under the conditions of Theorem~\ref{thuniversal-optimality}
least squares estimation combined with an optimal design can never be
improved by weighted least squares estimation.
In Section~\ref{general} several applications are presented.
In particular, we show that for a trigonometric system of regression
functions involving only cosinus terms with an arbitrary periodic
covariance kernel, the uniform distribution is universally optimal. We
also prove that the arcsine design is universally optimal
for the polynomial regression model with the logarithmic covariance kernel
and derive some universal optimality properties of the Beta distribution.
To our best knowledge these results provide the first explicit solution
of optimal design problems for regression models with correlated observations
which differ from the location scale model.

In Section~\ref{numerics} we provide an algorithm for computing
optimal designs
for any regression model with specified covariance function and investigate
the efficiency of the arcsine and uniform distribution in polynomial
regression models
with exponential correlation functions. Finally, Section~\ref{concl}
contains some conclusions
and technical details are given in the \hyperref[app]{Appendix}.

\section{Preliminaries} \label{prelim}

\subsection{The asymptotic covariance matrix}

Consider the linear regression model (\ref{lin-model}), where
$\ve(x)$ is a stochastic process with
%
%
\begin{equation}
\label{covkern} \mathrm{E}\ve(x)=0,\qquad \mathrm{E}\ve(x)\ve\bigl(x^\prime
\bigr)= K\bigl(x,x^\prime\bigr);\qquad x,x^\prime\in\mathcal{X}\subset
\mathbb{R}^d;
\end{equation}
the function $K(x,x^\prime)$ is called covariance kernel.
If $N$ observations,
say $y=(y_1,\ldots,y_N)^T$, are available
at experimental conditions $x_1,\ldots,x_N$ and the
covariance kernel is known, the vector of parameters can
be estimated by the weighted least squares method, that is,
$\hat\theta= (\mathbf{X}^T\bolds{\Sigma}^{-1}\mathbf
{X})^{-1}\mathbf
{X}^T\bolds{\Sigma}^{-1}y$,
where $\mathbf{X}=(f_i(x_j))^{i=1,\ldots,m}_{j=1,\ldots,N}$ and
$\bolds
{\Sigma}=(K(x_i,x_j))_{i,j=1,\ldots,N}$.
The variance--covariance matrix of this estimate is given by
\[
\operatorname{Var}(\hat\theta) =\bigl(\mathbf{X}^T\bolds{\Sigma
}^{-1}\mathbf{X}\bigr)^{-1}.
\]
If the correlation structure of the process is not known, one usually
uses the ordinary
least squares estimate $\tilde\theta= (\mathbf{X}^T\mathbf
{X})^{-1}\mathbf{X}^Ty$, which has the
covariance matrix
%
%
\begin{equation}
\label{eqlse} \operatorname{Var}(\tilde\theta) = \bigl(\mathbf{X}^T
\mathbf{X}\bigr)^{-1} \mathbf{X}^T\bolds{\Sigma} \mathbf{X}
\bigl(\mathbf{X}^T\mathbf{X}\bigr)^{-1}.
\end{equation}

An {exact} experimental design $\xi_N=\{x_1,\ldots,x_N\}$ is a
{collection} of
$N$ points in~$\mathcal{X}$, which defines the time points
or experimental conditions where observations are taken. Optimal
designs for weighted or ordinary least squares estimation
minimize a functional of the covariance matrix of
the weighted or ordinary least squares estimate, respectively, and
numerous optimality criteria have been proposed in the literature to
discriminate
between competing designs; see \citet{pukelsheim2006}.

Note that the weighted least squares
estimate can only be used if the correlation structure of the errors is known,
and its misspecification can lead to a considerable loss of efficiency.
At the same time, the ordinary least squares estimate does not employ
the structure of
the correlation. Obviously the ordinary least squares estimate
can be less efficient than the weighted least squares estimate, but in many
cases the loss of efficiency is small.
For example, consider the location scale model (\ref{const}) with a
stationary error process, the Gaussian correlation function $\rho
(t)=e^{-\lambda t^2}$
and the exact design $\xi=\{-1,-2/3,-1/3,1/3,2/3,1\}$.
{ Suppose that the guessed value of $\lambda$ equals $1$ while the true
value is $2$.
Then the variance of the weighted least squares estimate is $0.528$
computed as
\[
\bigl(\mathbf{X}^T\bolds{\Sigma}^{-1}_{\mathrm{guess}}
\mathbf{X}\bigr)^{-1} \mathbf{X}^T\bolds{
\Sigma}^{-1}_{\mathrm{guess}}\bolds{\Sigma}_{\mathrm{true}}\bolds{
\Sigma}^{-1}_{\mathrm{guess}} \mathbf{X} \bigl(\mathbf{X}^T
\bolds{\Sigma}^{-1}_{\mathrm{guess}}\mathbf{X}\bigr)^{-1},
\]
while the variance of the ordinary least squares estimate is $0.433$.
If the guessed value of $\lambda$ equals the true value, then
the variance of the weighted least squares estimate is $0.382$. }
A similar relation between the variances holds if
the location scale model and the Gaussian correlation function are replaced
by a polynomial model and a triangular or exponential correlation
function, respectively. For
a more detailed discussion concerning advantages of the ordinary least
squares against
the weighted least
squares estimate, see \citet{bickherz1979} and Section 5.1 in
\citet{naether1985a}.

Throughout this article we will concentrate on optimal
designs for the ordinary least squares estimate.
These designs also require the specification of the correlation
structure but a potential
loss by its misspecification in the stage of design construction is
typically much smaller
than the loss caused by the misspecification of the
correlation structure in the weighted least squares estimate.
Moreover, in this paper we will demonstrate that there are many
situations, where the combination of the ordinary least squares
estimate with the corresponding (universally) optimal design yields the
same covariance matrix as the weighted least squares estimate on the
basis of a (universally) optimal design for weighted least squares
estimation; see the discussions in Sections~\ref{general} and
\ref{concl}.

Because even in simple models the exact optimal designs are difficult
to find, most authors
usually use asymptotic arguments to determine efficient designs
for the estimation of the model parameters; see Sacks and Ylvisaker
(\citeyear{sackylv1966,sackylv1968}), \citet{bickherz1979}
or \citet{zhidetpep2010}.
Sacks and Ylvisaker (\citeyear{sackylv1966,sackylv1968}) and
N{\"a}ther [(\citeyear{naether1985a}), Chapter 4],
assumed that
the design points $\{x_{1},\ldots, x_{N}\}$
are generated by the quantiles of a distribution function, that is,
%
%
\begin{equation}
\label{despoint} x_{i}=a \bigl((i-1)/(N-1) \bigr),\qquad i=1,\ldots,N,
\end{equation}
where the function $a\dvtx[0,1]\to\mathcal{X}$ is the inverse of a
distribution function. If $\xi_N$ denotes a design with $N$ points and
corresponding\vspace*{1pt} quantile function $a(\cdot)$, the covariance
matrix of the least squares estimate $\tilde\theta=\tilde
\theta_{\xi_N}$ given in (\ref{eqlse}) can be written as
%
%
\begin{equation}
\label{eqvariance} \operatorname{Var}(\tilde\theta)=\bD(\xi_N) =
\bM^{-1}(\xi_N)\bB(\xi_N,\xi_N)
\bM^{-1}(\xi_N),
\end{equation}
where
%
%
\begin{eqnarray}
\label{inf} \bM(\xi_N)&=&\int_\mathcal{X} \bff(u)
\bff^T(u)\xi_N(\dd u),
\\
\label{infcorr} \bB(\xi_N,\xi_N)&=&\iint K(u,v)
\bff(u)\bff^T(v)\xi_N(\dd u)\xi_N(\dd v),
\end{eqnarray}
and $f(u)= (f_1(u),\ldots,f_m(u) )^T$ denotes the vector of
regression functions.
Following \citet{kiefer1974} we call any probability measure $\xi$ on
$\mathcal{X}$ (more precisely on an appropriate Borel field)
an approximate design or simply design.
The definition of the matrices $M(\xi)$ and
$B(\xi,\xi)$ can be extended to an arbitrary design $\xi$, provided that
the corresponding integrals exist.
The matrix
%
%
\begin{equation}
\label{eqDmatrix} \bD(\xi) = \bM^{-1}(\xi)\bB(\xi,\xi)
\bM^{-1}(\xi)
\end{equation}
is called the covariance matrix for the design $\xi$ and
can be defined for any probability measure $\xi$ supported on the
design space $\mathcal{X}$ such that the matrices $\bB(\xi,\xi)$ and
$\bM^{-1}(\xi)$ are well defined. This set will be denoted by $\Xi$.
An (approximate) optimal design
minimizes a functional of the covariance matrix $\bD(\xi)$ over the set
$\Xi$
and a universally optimal design $\xi^*$ (if it exists) minimizes the
matrix $\xi$ with respect to the Loewner ordering, that is,
\[
D\bigl(\xi^*\bigr) \leq D(\xi) \qquad\mbox{for all } \xi\in\Xi.
\]
Note that on the basis of this asymptotic analysis the kernel $K(u,v)$
has to be well defined for all $u,v \in\mathcal{X}$.
On the other hand, \citet{zhidetpep2010} extended the approach in
\citet
{bickherz1979}
and proposed an alternative approximation for the covariance matrix in
(\ref{eqlse}),
where the variance of the observations also depends on the sample size.
As a result they obtained an approximating matrix of the form (\ref{infcorr}),
where the kernel $K(u,v)$ in the matrix $B(\xi, \xi)$ may have
singularities at the diagonal.

Note that in general the function $\bD(\xi)$ is not
convex (with respect to the Loewner ordering) on the
space of all approximate designs. This implies that even if one
determines optimal designs by minimizing a convex
functional, say $\Phi$, of the matrix $D(\xi)$, the corresponding functional
$\xi\to\Phi(\bD(\xi))$ is generally not
convex on the
space of designs $\Xi$.
Consider, for example, the case $m=1$ where $\bD(\xi)$ is given by
%
%
\begin{equation}
\label{eqDDD} D(\xi)={ \biggl[ \int f^2(u)\xi(\dd u)
\biggr]^{-2}} {\iint K(u,v)f(u)f(v)\xi(\dd u)\xi(\dd v)},
\end{equation}
and it is obvious that this functional is not necessarily convex.
On the other hand, for the location scale model (\ref{const})
we have $m=1$, $f(x)=1$ for all $x$ and this expression reduces to
$
D(\xi)=\iint K(u,v)\xi(\dd u)\xi(\dd v).
$
In the stationary case $K(u,v)=\sigma^2\rho(u-v)$, where $\rho(\cdot)$
is a correlation function,
this functional is convex on the set of all probability measures on the
domain $\mathcal{X}$;
see Lemma 1 in \citet{zhidetpep2010} and Lemma 4.3 in \citet
{naether1985a}.
For this reason [namely the convexity of the functional $ D(\xi)$]
most of the literature discussing asymptotic optimal design problems
for least squares estimation in the presence of correlated observations
considers the location scale model, which corresponds to the estimation of
the mean of a stationary process; see, for example, \citet{bolnat1982},
N{\"a}ther (\citeyear{naether1985a,naether1985b}).

\subsection{Covariance kernels}
Consider the covariance kernels $K(u,v)$ that appeared in (\ref{covkern}).
An important case appears when the error process is stationary and the
covariance kernel
is of the form $K(u,v)= \sigma^2\rho(u-v)$,
where $\rho(0)=1$ and $\rho(\cdot)$ is called the correlation function.

Because in this paper we are interested in designs maximizing
functionals of the matrix $D(\xi)$
independently of the type of approximation which has been used to
derive it,
we will also consider singular kernels in the following discussion.
Moreover, we call $K(u,v)$ covariance kernel even if it has
singularities at the diagonal.

The covariance kernels with singularities at the diagonal can be used
as approximations to the
standard covariance kernels. They naturally appear as limits of
sequences of covariance kernels satisfying
$ K_N(u,v)=\sigma_N^2\rho_N(u-v),
$
where
$\rho_N(t) = \rho(a_N t)$, $\sigma_N^2 = a^\alpha_N \tau^2$,
$\tau>0$, $0<\alpha\leq1$, is a constant depending on
the asymptotic behavior of the function
$\rho(t)$ as $t \to\infty$,
and $\{ a_N\}_{N \in\mathbb{N}}$ denotes\vadjust{\goodbreak} a sequence of positive numbers
satisfying $ a_N \to\infty$ as $N \to\infty$.
Consider, for example, the correlation function $\rho
(t)=1/(1+|t|)^\alpha$ which is nonsingular. Then the sequence of functions
\[
\sigma_N^2\rho_N(t)=a_N^\alpha
\tau^2\frac{1}{(1+|a_Nt|)^\alpha
}=\tau^2\frac{1}{(1/a_N+|t|)^\alpha}
\]
converges to
$r_{\alpha}(t)=1/|t|^\alpha$ as $ N \to\infty$. For slightly different
types of approximation, see Examples~\ref{exam4} and~\ref{exam4a} below.

Let us summarize the assumptions regarding the covariance kernel.
First, we assume that $K$ is symmetric and continuous
at all points $(u,v) \in\mathcal{X} \times\mathcal{X}$ except
possibly at the diagonal points $(u,u)$.
We also assume that \mbox{$K(u,v)\neq0$}
for at least one pair $(u,v)$ with $u \neq v$.
Any covariance kernel $K(\cdot,\cdot)$ considered in this paper
is assumed to be positive definite in the following sense:
for any signed measure $\nu(\dd u)$ on $\mathcal{X}$, we have
%
%
\begin{equation}
\label{eqposdef} \iint K(u,v) \nu(\dd u) \nu(\dd v) \geq0.
\end{equation}
If the kernel $K(u,v)$ has singularities at the diagonal then the
assumptions we make are as follows. We assume that
$K(u,v)= r(u-v)$,
where $r(\cdot)$ is a function on $\Real\setminus\{0\} $ with $ 0
\leq r(t) < \infty$ for all $t\neq0$ and
$r(0)=+\infty$. We also assume that
there exists a
monotonously increasing
sequence $\{\sigma^2_N \rho_N(t)\}_{N \in\mathbb{N}}$
of covariance functions such that
$0\leq\sigma^2_N \rho_N(t) \leq r(t)$ for all $t$ and all
$N=1,2,\ldots
$ and
$r(t) = \lim_{N \to\infty} \sigma^2_N \rho_N(t)$. Theorem 5 in
\citet
{zhidetpep2010} then guarantees that
for this kernel we also have the property of positive definiteness~(\ref{eqposdef}).


%

\subsection{The set of admissible design points}

Consider the vector-function $f(x) = (f_1(x), \ldots, f_m(x))^T$ used
in the definition of the regression model (\ref{lin-model}). Define the
sets $\mathcal{X}_0=\{x\in\mathcal{X}\dvtx f(x)=0\}$
and $\mathcal{X}_1=\mathcal{X}\setminus\mathcal{X}_0=\{x\in
\mathcal
{X}\dvtx\break f(x) \neq0\}$ and assume that designs $\xi_0$ and $\xi_1$ are
concentrated on $\mathcal{X}_0$ and $\mathcal{X}_1$
correspondingly.
Consider the design $\xi_\alpha=\alpha\xi_0+(1-\alpha)\xi_1$ with
$0\le
\alpha<1$;
note that if the design $\xi_\alpha$ is concentrated on the set
$\mathcal{X}_0$ only (corresponding to the case $\alpha=1$),
then the construction of estimates is not possible.
We have
\[
M(\xi_\alpha)=\int f(x)f^T(x)\xi_\alpha(\dd x)=(1-
\alpha)M(\xi_1),\qquad M^{-1}(\xi_\alpha)=
\frac{1}{1-\alpha}M^{-1}(\xi_1)
\]
and
\[
B(\xi_\alpha,\xi_\alpha)=\iint K(x,u)f(x)f^T(u)
\xi_\alpha(\dd x)\xi_\alpha(\dd u)=(1-\alpha)^2B(
\xi_1,\xi_1).
\]
Therefore, for all $0\le\alpha<1$ we have
\[
D(\xi_\alpha)=M^{-1}(\xi_\alpha)B(\xi_\alpha,
\xi_\alpha)M^{-1}(\xi_\alpha) =M^{-1}(
\xi_1)B(\xi_1,\xi_1)M^{-1}(
\xi_1)= D(\xi_1).
\]
Consequently,\vspace*{1pt} observations taken at points from the set
$\mathcal{X}_0$ do not change the estimate $\hat\theta$ and its
covariance matrix. If we use the convention\vadjust{\goodbreak} $0\cdot\infty=0$, it
follows that $\iint K(x,u)f(x)f^T(u) \xi_0(\dd x) \xi_0(\dd u)=0$, and
this statement is also true for the covariance kernels $K(x,u)$ with
singularity at $x=u$.

Summarizing this discussion, we assume throughout this paper that
\mbox{$f(x)\! \neq\!0$} for all $x\in\mathcal{X}$.


\section{Characterizations of optimal designs}
\label{equiv}

\subsection{General optimality criteria}\label{sec3.1}
Recall the definition of the information matrix in (\ref{inf}) and define
\[
\bB(\xi,\nu)=\int_\mathcal{X} \int_\mathcal{X}
K(u,v)\bff(u)\bff^T(v)\xi(\dd u)\nu(\dd v),
\]
where $\xi$ and $\nu\in\Xi$ are two arbitrary designs,
and $K(u,v)$ is a covariance kernel.

According to the discussion in the previous paragraph,
the asymptotic covariance matrix of the least squares estimator $\hat
{\theta}$ is proportional to
the matrix $D(\xi)$ defined in (\ref{eqvariance}).
Let $ \Phi(\cdot)$ be a monotone, real valued functional defined on the
space of symmetric $m\times m$ matrices
where the monotonicity of $\Phi(\cdot)$ means that $A\geq B$ implies
$\Phi(A) \geq\Phi(B)$.
Then the optimal design $\xi^*$ minimizes the function
%
%
\begin{equation}
\label{critgen} \Phi\bigl(\bD(\xi)\bigr)
\end{equation}
on the space $\Xi$ of all designs.
In addition to monotonicity, we shall also assume differentiability of
the functional $\Phi(\cdot)$; that is,
the existence of the matrix of derivatives
\[
\bC=\frac{\partial\Phi(D)}{\partial D}= \biggl( \frac{\partial
\Phi
(D)}{\partial D_{ij}} \biggr)_{i,j=1,\ldots,m},
\]
where $D$ is any symmetric nonnegative definite matrix of size $m
\times m$.
The following lemma is crucial in the proof of the optimality theorem below.

%
\begin{lemma}
Let $\xi$ and $\nu$ be two designs and $\Phi$ be a differentiable functional.
Set $\xi_\alpha=(1-\alpha)\xi+\alpha\nu$ and assume that the matrices
$\bM(\xi)$ and $\bB(\xi,\xi)$ are nonsingular.
Then the directional derivative of $\Phi$ at the design $\xi$ in the
direction of $\nu- \xi$ is given by
\[
\frac{\partial\Phi(\bD(\xi_\alpha))}{\partial\alpha}\bigg|_{\alpha=0}= 2
\bigl[ {\mathbf b}(\nu,\xi)- \bvarphi(\nu,
\xi)\bigr],
\]
where
\begin{eqnarray*}
\bvarphi(\nu,\xi)&=&\tr\bigl({\bM(\nu) \bD(\xi) \bC(\xi) \bM^{
-1}(\xi)
}\bigr),
\\
{\mathbf b}(\nu,\xi)&=&\tr\bigl(\bM^{ -1}(\xi) \bC(\xi)
\bM^{ -1}(\xi) \bB(\xi,\nu)\bigr)
\end{eqnarray*}
and
\[
\bC(\xi)=\frac{\partial\Phi(D)}{\partial D}\bigg|_{D=\bD(\xi)}.\vadjust{\goodbreak}
\]
\end{lemma}
\begin{pf}
Straightforward calculation shows that
\[
\frac{\partial}{\partial\alpha}\bM^{ -1}(\xi_\alpha) \bigg|_{\alpha=0}
=\bM^{ -1}(\xi)- \bM^{ -1}(
\xi)\bM(\nu)\bM^{ -1}(\xi)
\]
and
%
\[
\frac{\partial}{\partial\alpha}\bB(\xi_\alpha,\xi_\alpha) \bigg|_{\alpha=0}
=\bB(\xi,\nu)+\bB(\nu,\xi)-2\bB(
\xi,\xi).
\]
Using the formula for the derivative of a product and the two formulas above,
we obtain
\[
\frac{\partial}{\partial\alpha}\bD(\xi_\alpha) \bigg|_{\alpha=0} 
= - 2 \bM^{ -1}(\xi)\bM(\nu)\bD(\xi) +2 \bM^{ -1}(
\xi) \bB(\xi,\nu) \bM^{ -1}(\xi).
\]
Note that the matrices $\bM(\xi_{\alpha})$ and $B(\xi_\alpha, \xi_\alpha
)$ are nonsingular
for small nonnegative $\alpha$ (i.e., for all $\alpha\in[0,\alpha_0)$
where $\alpha_0$ is a small positive number)
which follows from the nondegeneracy of $\bM(\xi)$ and $B(\xi, \xi)$
and the continuity of $\bM(\xi_{\alpha})$ and $B(\xi_\alpha, \xi_\alpha
)$ with respect to $\alpha$.

Using the above formula and the fact that
$\tr(H(A+A^T))=2\tr(HA)$ for any $m \times m$ matrix $A$ and any $m
\times m$ symmetric matrix $H$, we obtain
\[
\frac{\partial\Phi(\bD(\xi_\alpha))}{\partial\alpha}\bigg|_{\alpha=0} = \tr
\biggl( C(\xi) \,\frac{\partial}{\partial\alpha}\bD(
\xi_\alpha) \biggr) \bigg|_{\alpha=0} 
=
2\bigl[{\mathbf b}(\nu,\xi) - \bvarphi(\nu,\xi)\bigr].
\]
\upqed\end{pf}

Note that the functions ${\mathbf b}(\nu,\xi)$ and $ \bvarphi(\nu,\xi)$
can be represented as
\[
{\mathbf b}(\nu,\xi)= \int b(x,\xi) \nu(\dd x),\qquad \bvarphi(\nu,\xi)= \int
\varphi(x,\xi) \nu(\dd x),
\]
where
%
%
\begin{eqnarray}
\label{eqphi-x-xi} \varphi(x,\xi)&=& \bvarphi(\xi_x,\xi)= {
\bff^T(x) \bD(\xi) \bC(\xi) \bM^{ -1}(\xi)\bff(x)},
\\
\label{eqb-x-xi} b(x,\xi)&=& {\mathbf b} (\xi_x, \xi) = \tr\bigl(C(
\xi)M^{-1}(\xi)B(\xi,\xi_x)M^{-1}(\xi)\bigr),
\end{eqnarray}
and $\xi_x$ is the probability measure concentrated at a point $x$.

%
\begin{lemma}\label{lemma3.2}
For any design $\xi$ such that the matrices $\bM(\xi)$ and
$\bB(\xi,\xi)$ are nonsingular we have
%
%
\begin{equation}
\label{eqint-phi=tr} \int\varphi(x,\xi)\xi(\dd x)=\int b(x,\xi)\xi(\dd
x)=\tr\bD(
\xi) \bC(\xi),
\end{equation}
where the functions $\varphi(x,\xi)$ and $b(x,\xi)$ are defined in
(\ref{eqphi-x-xi})
and (\ref{eqb-x-xi}), respectively.
\end{lemma}
\begin{pf}
Straightforward calculation shows that
\[
\int\varphi(x,\xi)\xi(\dd x) 
= \tr\biggl(\bD(\xi) \bC(
\xi) \bM^{ -1}(\xi) \int f(x)f^T(x) \xi(\dd x)\biggr) = \tr
\bigl(\bD(\xi) \bC(\xi)\bigr).\vadjust{\goodbreak}
\]
We also have
\begin{eqnarray*}
\int\bB(\xi,\xi_x)\xi(\dd x) &=& \int\biggl[\iint K(u,v)\bff(u)
\bff^T(v)\xi(\dd u)\xi_x(\dd v) \biggr] \xi(\dd x)
\\[-2pt]
&=& \int\biggl[\int K(u,x)\bff(u)\bff^T(x)\xi(\dd u) \biggr] \xi(
\dd x)=\bB(\xi,\xi),
\end{eqnarray*}
which implies
\begin{eqnarray*}
\int b(x,\xi)\xi(\dd x) &=& \tr\biggl(\bM^{ -1}(\xi) \bC(\xi)
\bM^{
-1}(\xi) \int\bB(\xi,\xi_x)\xi(\dd x)\biggr)
\\
&=& \tr\bigl(\bD(\xi) \bC(\xi)\bigr).
\end{eqnarray*}
\upqed\end{pf}

The first main result of this section provides a necessary condition
for the optimality of a given design.

%
\begin{theorem}\label{thphi-opt-cond-pdim}
Let $\xi^*$ be any design minimizing the functional $\Phi(\bD(\xi))$.
Then the inequality
%
%
\begin{equation}
\label{eqphi-ekv-multi} \varphi\bigl(x,\xi^*\bigr)\le b\bigl(x,\xi
^*\bigr)
\end{equation}
holds for all $x \in\mathcal{X}$, where the functions $\varphi(x,\xi)$
and $b(x,\xi)$ are defined in (\ref{eqphi-x-xi})
and (\ref{eqb-x-xi}), respectively.
Moreover, there is equality in (\ref{eqphi-ekv-multi}) for $\xi
^*$-almost all $x$, that is,
$\xi^*(\mathcal{A})=0$ where
\[
\mathcal{A} = \mathcal{A}\bigl(\xi^*\bigr) = \bigl\{ x \in\mathcal{X}
\mid
\varphi\bigl(x,\xi^*\bigr) < b\bigl(x,\xi^*\bigr) \bigr\}
\]
is the set of $x \in\mathcal{X} $
such that the inequality (\ref{eqphi-ekv-multi}) is strict.
\end{theorem}
\begin{pf}
Consider any design $\xi^*$ minimizing the functional $\Phi(\bD(\xi))$.
The necessary condition for an element to be a minimizer of a
differentiable functional states
that the directional derivative from this element in any direction is
nonnegative.
In the case of the design $\xi^*$ and the functional $\Phi(\bD(\xi))$
this yields for any design $\nu$
\[
\frac{\partial\Phi(\bD(\xi_\alpha))}{\partial\alpha}\bigg|_{\alpha
=0}\ge0,
\]
where $\xi_\alpha=(1-\alpha)\xi^*+\alpha\nu$.
Inequality (\ref{eqphi-ekv-multi}) follows now from
Lemma 1.
The assumption that inequality (\ref{eqphi-ekv-multi}) is strict
for all $x \in\mathcal{A}$ with \mbox{$\xi^*(\mathcal{A})>0$}
is in contradiction with identity (\ref{eqint-phi=tr}).
\end{pf}

%
\begin{remark}
In the classical theory of optimal design, convex optimality criteria
are almost always considered. However, in at least one paper, namely
\citet{torsney1986}, an optimality theorem for a rather general
nonconvex optimality criteria was established and used (in the case of
noncorrelated observations).
\end{remark}

\subsection{An alternative representation of
the necessary condition of optimality}

For a given design $\xi\in\Xi$, introduce the vector-valued function
%
%
\begin{equation}
\label{eqintKfxi=Lf+g} g(x)= \int K (x,u) f(u) \xi(\dd u)- \Lambda
f(x),\qquad
x \in{
\cal X},
\end{equation}
where $\Lambda=B(\xi,\xi) M^{-1}(\xi) $.
This function satisfies the equality
%
%
\begin{equation}
\label{eqcond-orth-g} \int g(x)f^T(x)\xi(\dd x)=0.
\end{equation}
Additionally, as the vector of regression functions $f(\cdot)$ is
continuous on ${\cal X}$, the function $g(\cdot) $ is continuous too.

Note that $f_1,\ldots,f_m\in L_2(\mathcal{X},\xi)$ where
\[
L_2(\mathcal{X},\xi)= \biggl\{ h\dvtx\mathcal{X} \to\mathbb{R} \Big|
\int h^2(x)\xi(\dd x)<\infty\biggr\}.
\]
Formula (\ref{eqintKfxi=Lf+g}) implies that
$g(x)$ is the residual obtained after component-wise projection of the
vector-valued function $\int K(x,u)f(u) \xi(\dd u)$ onto the subspace
$\operatorname{span}\{f_1,\ldots,f_m\}\subset L_2(\mathcal{X},\xi)$.

Using (\ref{eqintKfxi=Lf+g}), (\ref{eqcond-orth-g}) and the symmetry
of the matrix $B(\xi,\xi)$ we obtain
\begin{eqnarray*}
B(\xi,\xi)&=&\iint K(x,u)f(u)\xi(\dd u)f^T(x)\xi(\dd x)
\\
&=&\int\Lambda f(x)f^T(x)\xi(\dd x)+\int g(x)f^T(x)\xi(
\dd x)
\\
&=&\Lambda M(\xi)=M(\xi)\Lambda^T,
\end{eqnarray*}
which gives for the matrix $D$ in (\ref{eqDmatrix}),
\[
\bD(\xi)= \bM^{ -1}(\xi) \bB(\xi,\xi)\bM^{ -1}(\xi)=
\bM^{
-1}(\xi) \Lambda= \Lambda^T \bM^{ -1}(\xi).
\]
For the function (\ref{eqphi-x-xi}), we obtain
\begin{eqnarray*}
\varphi(x,\xi)&=&{\bff^T(x) \bD(\xi) \bC(\xi) \bM^{ -1}(\xi
)\bff(x)}
\\
&=&{\bff^T(x) \Lambda^T\bM^{ -1}(\xi) \bC(\xi)
\bM^{ -1}(\xi)\bff(x)}
\\
&=&{\bff^T(x) \bM^{ -1}(\xi) \bC(\xi) \bM^{ -1}(
\xi)\Lambda\bff(x)}.
\end{eqnarray*}
We also have
\[
B(\xi,\xi_x)=\int K(x,u)f(u)\xi(\dd u)f^T(x)=\Lambda
f(x)f^T(x)+g(x)f^T(x),
\]
which gives for function (\ref{eqb-x-xi})
\begin{eqnarray*}
b(x,\xi)&=&\tr\bigl(C(\xi)M^{-1}(\xi)B(\xi,\xi_x)M^{-1}(
\xi)\bigr)
\\
&=&\varphi(x,\xi)+r(x,\xi),
\end{eqnarray*}
where the function $r$ is defined by
\[
r(x,\xi)=f^T(x)M^{-1}(\xi)C(\xi)M^{-1}(\xi)g(x).
\]
The following result is now an obvious corollary of Theorem~\ref{thphi-opt-cond-pdim}.

%
\begin{corollary}\label{cor3.1}
If a design $\xi$ is optimal, then
$r(x,\xi)\ge0$ for all $x \in{\cal X}$ and $r(x,\xi)=0$ for all $x$
in the support of the measure $\xi$.
\end{corollary}

\subsection{$D$-optimality}

For the $D$-optimality there exists an analogue of the
celebrated ``equivalence theorem'' of \citet{kiewol1960}, which
characterizes optimal designs minimizing the $D$-optimality
criterion\break
$\Phi(\bD(\xi))= \ln\det(\bD(\xi))$.
%
%
\begin{theorem}\label{theo3.2}
Let $\xi^*$ be any $D$-optimal design. Then
for all $x \in\mathcal{X}$ we have
%
%
\begin{equation}
\label{eqekv-multi} d\bigl(x,\xi^*\bigr) \le b\bigl(x,\xi^*\bigr),
\end{equation}
where the functions $d$ and $b$ are defined by
$d(x,\xi)={f^T(x) \bM^{ -1}(\xi) f(x)}$ and
%
%
\begin{eqnarray}
\label{eqbdopt} b(x,\xi)&=&\tr\bigl( \bB^{-1}(\xi,\xi)\bB(\xi,
\xi_x) \bigr)\nonumber\\[-8pt]\\[-8pt]
&=& \bff^T(x)\bB^{-1}(\xi,\xi)\int
K(u,x)\bff(u)\xi(\dd u),\nonumber
\end{eqnarray}
respectively. Moreover, there is equality in (\ref{eqekv-multi}) for
$\xi^*$-almost all $x$.
\end{theorem}
\begin{pf}
In the case of the $D$-optimality criterion $\Phi(\bD(\xi))= \ln
\det
(\bD(\xi))$, we have
$
\bC(\xi)=\bD^{-1}(\xi),
$
which gives
\[
\varphi(x,\xi) = {\bff^T (x) \bD(\xi) \bD^{-1}(\xi)
\bM^{ -1} (\xi)\bff(x)} = d(x,\xi).
\]
Similarly, we simplify an expression for $b(x,\xi)$. Reference to
Theorem~\ref{thphi-opt-cond-pdim} completes the proof.
\end{pf}

Note that the function $r(x,\xi)$ for the $D$-criterion is given by
\[
r(x,\xi)=f^T(x)B^{-1}(\xi,\xi)g(x)
\]
and, consequently, the necessary condition of the $D$-optimality can be
written as $f^T(x)B^{-1}(\xi,\xi)g(x)\ge0$ for all $x\in\mathcal{X}$.

The following statement illustrates a remarkable similarity between
$D$-optimal design problems
in the cases of correlated and noncorrelated observations.
The proof easily follows from Lemma~\ref{lemma3.2} and Theorem~\ref{theo3.2}.
%
%
\begin{corollary}
For any design $\xi$ such that the matrices $\bM(\xi)$ and
$\bB(\xi,\xi)$ are nonsingular we have
\[
\int d(x,\xi)\xi(\dd x)=\int b(x,\xi)\xi(\dd x)=m,\vadjust{\goodbreak}
\]
where $b(x,\xi)$ is defined in (\ref{eqbdopt})
and $m$ is the number of parameters in the regression model (\ref{lin-model}).
\end{corollary}

%
\begin{example} \label{exam1}
Consider the quadratic regression model
$
y(x) = \theta_1 + \theta_2x + \theta_3 x^2 + \varepsilon(x)
$
with design space $\mathcal{X} = [-1,1]$.
In Figure~\ref{figdxbx} we plot
functions $b(x,\xi)$ and $d(x,\xi)$ for the covariance kernels
%
%
\begin{figure}

\includegraphics{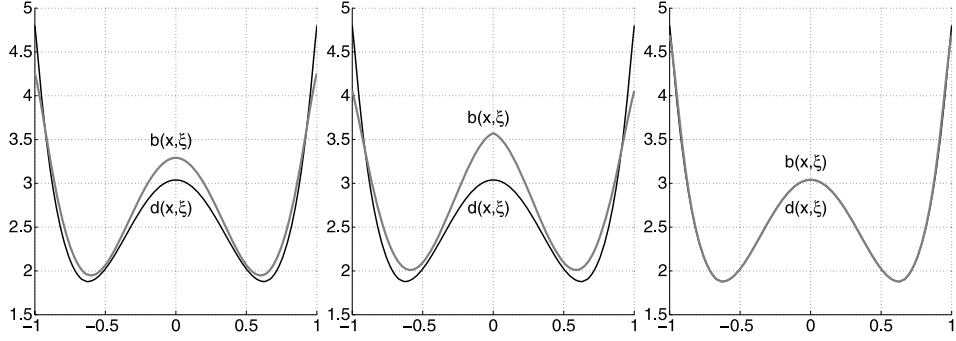}

\caption{The functions $b(x,\xi)$ and $d(x,\xi)$ for the regression
model (\protect\ref{lin-model}) with $f(x)=(1,x,x^2)^T$ and the covariance
kernels $K(u,v)=e^{-|u-v|}$ (left), $K(u,v)=\max(0,1-|u-v|)$ (middle)
and $K(u,v)=-\log(u-v)^2$ (right) and the arcsine design $\xi_{a}$.}
\label{figdxbx}
\end{figure}
$K(u,v) = e^{-|u-v|}$, $K(u,v) = \max\{0,1-|u-v| \}$ and $K(u,v) = -
\log(u-v)^2$,
where the design is the arcsine distribution 
with density
%
%
\begin{equation}
\label{arcs} p(x)=1/\bigl(\pi\sqrt{1-x^2}\bigr),\qquad x \in(-1,1).
\end{equation}
Throughout this paper this design will be called ``arcsine design'' and
denoted by~$\xi_a$.
By the definition,
the function $d(x,\xi)$ is the same for different covariance kernels, but
the function $b(x,\xi)$ depends on the choice of the kernel. From the
left and middle panel we see that the arcsine design
does not satisfy the necessary condition of Theorem \ref
{thphi-opt-cond-pdim} for
the kernels $K(u,v)=e^{-|u-v|}$ and $\max\{ 0,1-|u-v|\}$ and is
therefore not $D$-optimal for the quadratic regression model.
On the other hand, for the logarithmic kernel $K(u,v)=- \log(u-v)^2$
the necessary condition is satisfied, and the arcsine design $\xi_a$
is a candidate for the $D$-optimal design.
We will show in Theorem~\ref{THmulti-genarcsine45} that the design
$\xi_a$
is universally optimal and as a consequence optimal with respect to
a broad class of criteria including the $D$-optimality criterion.
\end{example}


\subsection{$c$-optimality}

For the $c$-optimality criterion $\Phi(D(\xi))=c^TD(\xi)c$,
we have $C(\xi)=cc^T$.
Consequently,
\[
\varphi(x,\xi)=\bff^T(x) \bM^{ -1}(\xi) cc^T
\bM^{ -1}(\xi)\Lambda\bff(x) =c^T \bM^{ -1}(\xi)
\Lambda\bff(x)\bff^T(x) \bM^{ -1}(\xi) c
\]
and
\[
r(x,\xi) = b(x,\xi)-\varphi(x,\xi) = f^T(x)M^{-1}(\xi)c
c^TM^{-1}(\xi)g(x).
\]
Therefore, the necessary condition for $c$-optimality simplifies to
%
%
\begin{equation}
\label{eqc-opt} f^T(x)M^{-1}(\xi)c c^TM^{-1}(
\xi)g(x)\ge0 \qquad\mbox{for all $x\in\mathcal{X}$.}
\end{equation}

%
\begin{example} \label{exam2}
Consider again the quadratic regression model $ y(x) = \theta_1 +
\theta_2x + \theta_3 x^2 + \varepsilon(x) $ with design space
$\mathcal{X} = [-1,1]$. Assume the triangular correlation function
$\rho(x)=\max\{0,1-|x|\}$.

Let $\xi=\{-1,0,1; 1/3,1/3,1/3\}$ be the design assigning weights
$1/3$ to the points $-1, 0$ and 1.
For this design, we have the matrices $M(\xi)$ and $D(\xi)$
\[
M(\xi)=\pmatrix{ 1 & 0 & 2/3
\cr
0 & 2/3 & 0
\cr
2/3 & 0 & 2/3},\qquad D(\xi)=
\pmatrix{ 1&0&-1
\cr
0&1/2&0
\cr
-1&0&3/2},
\]
and the matrix $\Lambda$ and the vector $g$ are given by
\[
\Lambda=\operatorname{diag}(1/3,1/3,1/3),\qquad g(x)=\bigl(1/3,x/3,|x|/3\bigr)^T.
\]

%
%
\begin{figure}

\includegraphics{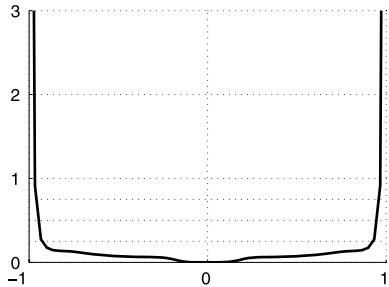}

\caption{The $c$-optimal design for the quadratic model and the
triangular correlation function, where $c=(1,0,0)^T$.}
\label{fige1optdes}
\end{figure}

%
%
\begin{figure}

\includegraphics{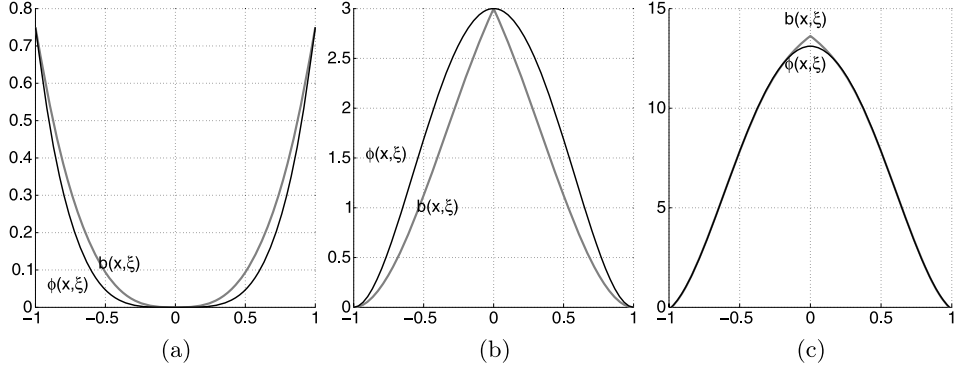}

\caption{The functions $b(x,\xi)$ and $\phi(x,\xi)$ for the
$c$-optimality criterion. \textup{(a)}: $c=(1,0,1)^T$, design $\xi=\{-1,0,1;
1/3,1/3,1/3\}$; \textup{(b)}: $c=(1,0,0)^T$, design
$\xi=\{-1,0,1;1/3,1/3,1/3\}$; \textup{(c)}:~$c=(1,0,0)^T$, design is displayed
in Figure \protect\ref{fige1optdes}.} \label{figcopt}
\end{figure}

If $c=(0,1,0)^T$, then $r(x,\xi)=0$ for all $x \in[-1,1]$ and thus the
design $\xi$ satisfies the necessary condition for $c$-optimality in
(\ref{eqc-opt}).
If $c=(1,0,1)^T$, then $r(x,\xi)=\frac34 |x|^3(1-|x|)\geq0$ for all $x
\in[-1,1]$ and
the design $\xi$ also satisfies (\ref{eqc-opt}). The corresponding
functions $b$ and $\varphi$ are displayed in the left and middle panels
of Figure~\ref{figcopt}.
Numerical analysis shows that for both vectors this design is in fact
$c$-optimal. However, it is not optimal for
any $c$-optimality criteria.
For example, if $c=(1,0,0)^T$, then $r(x,\xi)=
-3 x ( 1-|x| ) ( 1-{x}^{2} ) \leq0$ for all $x
\in[-1,1]$, showing that
the design is not $c$-optimal; see the middle panel of Figure~\ref
{figcopt}. For this case, the density function of the $c$-optimal design
is
displayed in Figure~\ref{fige1optdes}. The corresponding
functions $b$ and $\varphi$ are shown in the right panel of Figure
\ref{figcopt}.
\end{example}

\subsection{Universal optimality}

In this section we consider the matrix $D(\xi)$ defined in (\ref
{eqDmatrix}) as the matrix optimality criterion which we are going to minimize
on the set $\Xi$
of all designs, such that the
matrices $\bB(\xi,\xi)$ and $\bM^{-1}(\xi)$ [and therefore the matrix
$D(\xi)$] are well defined.
Recall that a design $\xi^*$ is universally optimal if
$
\bD(\xi^*)\le\bD(\xi)
$
in the sense of the Loewner ordering for any design $\xi\in\Xi$.
Note that a design $\xi^*$ is universally optimal
if and only if $\xi^*$ is $c$-optimal for any vector $c \in\mathbb
{R}^m \setminus\{ 0 \}$;
that is,
$c^T D(\xi^*) c \le c^T D(\xi) c$ for any $\xi\in{\Xi}$ and any $c
\in\mathbb{R}^m $.

%
\begin{theorem}
\label{thuniversal-optimality}
Consider the regression model (\ref{lin-model}) with a covariance
kernel~$K$, a design $\xi\in\Xi$ and the corresponding
the vector-function $g(\cdot)$ defined in~(\ref{eqintKfxi=Lf+g}).

\begin{longlist}[(a)]
\item[(a)]
If $g(x)=0$ for all $x \in\mathcal X$, then the design
$\xi$ is universally optimal;

\item[(b)] If
the design
$\xi$ is universally optimal, then the function $g(\cdot)$
can be represented in the form $g(x)=\gamma(x) f(x)$, where $\gamma
(x)$ is a nonnegative function defined on $\mathcal{X}$ such that
$\gamma(x)=0$ for all $x$ in the support of the design~$\xi$.
\end{longlist}
\end{theorem}


In the proof of Theorem~\ref{thuniversal-optimality} we shall need the
following two auxiliary results which will be proved in the \hyperref
[app]{Appendix}.

%
\begin{lemma}
\label{lem-convex}
Let $c\in\mathbb{R}^m$, and $\mathcal{M}$ be the set of
all signed vector measures supported on $\mathcal{X}$.
Then the functional $\Phi_c\dvtx \mathcal{M} \to\mathbb{R}_+$ defined by
%
%
\begin{equation}
\label{auxfct} \Phi_c(\mu) = c^T \iint K(x,u) \mu(\dd
x) \mu^T (\dd u) c
\end{equation}
is convex.
\end{lemma}

%
\begin{lemma}
\label{lemcabc<0}
Let $m>1$ and $a,b\in\mathbb{R}^m$ be two linearly independent vectors.
Then there exists a vector $c\in\mathbb{R}^m$ such that
$S_c=c^Tab^Tc<0.$\vadjust{\goodbreak}
\end{lemma}

\begin{pf*}{Proof of Theorem~\ref{thuniversal-optimality}}
Consider the regression model $y(x)= f^T(x) \theta+ \varepsilon(x)$,
where the full trajectory $\{ y(x) | x \in\mathcal{X} \}$ can be observed.
Let
$
\hat\theta(\mu) = \int y(x) \mu(\dd x)
$
be a general linear unbiased estimate of the parameter $\theta$, where
$\mu=(\mu_1,\ldots,\mu_m)^T$ is a vector of signed measures. For
example, the least squares estimate for a design $\xi$ in this model is
obtained as $\hat\theta(\mu_\xi)$, where $\mu_\xi(\dd x)=M^{-1}(\xi
)f(x)\xi(\dd x)$.
The condition of unbiasedness of the estimate $\hat\theta(\mu)$
means that
\[
\theta= \mathrm{E}\bigl[\hat\theta(\mu)\bigr] = \mathrm{E} \biggl[\int
\mu(\dd
x) y(x) \biggr] = \int\mu(\dd x) f^T(x) \theta
\]
for all $\theta\in\mathbb{R}^m$, which is equivalent to the condition
%
%
\begin{equation}
\label{rest} \int\mu(\dd x) f^T(x) = \int f(x) \mu^T (\dd
x) = I_m,
\end{equation}
where $I_m$ denotes the $m \times m$ identity matrix.
In the following discussion we define $\mathcal{M}_0$ as a subset of
$\mathcal{M}$ containing the signed
measures which
satisfy condition~(\ref{rest}). Note that both sets, $\mathcal{M}$ and
$\mathcal{M}_0$, are convex.

For a given vector $c \in\mathbb{R}^m$, the variance of the estimate
$c^T \hat\theta(\mu)$ is given by
\begin{eqnarray*}
\operatorname{Var} \bigl(c^T \hat\theta(\mu)\bigr) &=&
c^T \iint\mathrm{E}\bigl[\varepsilon(x)\varepsilon(u)\bigr] \mu(
\dd x) \mu^T (\dd u)c
\\
&=& c^T \iint K(x,u) \mu(\dd x) \mu^T (\dd u) c =
\Phi_c(\mu),
\end{eqnarray*}
and a minimizer of this expression with respect to $\mu\in\mathcal{M}_0$
determines the best linear unbiased estimate for $c^T \theta$ and the
corresponding $c$-optimal design simultaneously.

Note that the sets $\mathcal{M}$ and $\mathcal{M}_0 $ are convex and
in view of Lemma~\ref{lem-convex} the functional
$\Phi_c(\mu)$ defied in (\ref{auxfct}) is convex on $\mathcal{M}$.
Similar arguments as given in Section~\ref{sec3.1} show that
the directional derivative of $\Phi_c$ at $\mu^*$ in the direction of
$\nu- \mu^*$ is given by
\begin{eqnarray*}
&&
\frac{\partial}{\partial\alpha} \Phi_c (\mu_\alpha) \bigg|_{\alpha=
0} \\
&&\qquad=
\frac{\partial}{\partial\alpha} \Phi_c \bigl((1 - \alpha) \mu^* +
\alpha\nu\bigr)
\bigg|_{\alpha= 0}
\\
&&\qquad= 2c^T \biggl[ \iint K(x,u) \mu^* (\dd x) \nu^T (\dd
u) - \iint K(x,u) \mu^* (\dd x) \mu^{*T}(\dd u) \biggr] c.
\end{eqnarray*}
Because $\Phi_c$ is convex, the optimality of $\mu^*$ in the set
$\mathcal{M}_0 $ is equivalent to
the condition $\frac{\partial}{\partial\alpha} \Phi_c (\mu_\alpha)
|_{\alpha= 0} \geq0$ for all $\nu\in\mathcal{M}_0$.
Therefore, the signed measure $\mu^* \in\mathcal{M}_0 $ minimizes the
functional $\Phi_c(\mu)$ if and only if the inequality
%
%
\begin{equation}
\label{check} \Phi_c \bigl(\mu^*,\nu\bigr) \geq\Phi_c
\bigl(\mu^*\bigr)
\end{equation}
holds for all $\nu\in\mathcal{M}_0$, where
\[
\Phi_c \bigl(\mu^*,\nu\bigr):= c^T \iint K (x,u)
\mu^* (\dd x) \nu^T (\dd u) c.
\]

Let us prove part (a) of Theorem~\ref{thuniversal-optimality}.
Consider a design $\xi\in\Xi$ such that \mbox{$g(x)=0$} for all $x \in
{\cal X}$
and define
the vector-valued measure $\mu_0(\dd x)= M^{-1}(\xi) \*f(x)\xi(\dd x)$.
It follows for all $\nu\in\mathcal{M}_0$,
\begin{eqnarray*}
\Phi_c (\mu_0,\nu) &=& c^T \iint K (x,u)
M^{-1}(\xi) f(x) \xi(\dd x) \nu^T (\dd u) c
\\
&=& c^T M^{-1}(\xi) \Lambda\int f(u) \nu^T (
\dd u) c = c^T M^{-1}\bigl(\xi^*\bigr) \Lambda c,
\end{eqnarray*}
where we used (\ref{rest}) for the measure $\nu$ in the last identity.
On the other hand, $\mu_0 \in\mathcal{M} $ also satisfies
(\ref{rest}),
and we obtain once more using identity (\ref{eqintKfxi=Lf+g}) with
$g(x) \equiv0$,
\begin{eqnarray*}
\Phi_c (\mu_0) &=& c^T \iint K(x,u)
\mu_0 (\dd x) \mu_0^T (\dd u) c
\\
&= & c^T \int\biggl[ \int K(x,u) M^{-1}(\xi) f(x) \xi(\dd
x) \biggr] \mu_0^T (\dd u) c
\\
&= & c^T M^{-1}(\xi) \Lambda\int f(u)
\mu_0^T (\dd u) c = c^T M^{-1}(\xi)
\Lambda c.
\end{eqnarray*}
This yields that for $\mu^*=\mu_0$ we have equality in (\ref{check})
for all $\nu\in\mathcal{M}_0$, which shows that the vector-valued
measure $\mu_0(\dd x)=
M^{-1}(\xi) f(x)\xi(\dd x)$ minimizes the function $\Phi_c$
for any $c \neq0$ over the set $\mathcal{M}_0$ of signed vector-valued
measures satisfying (\ref{rest}).

Now we return to the minimization of the function $D(\eta)$ in the
class of all designs $\eta\in\Xi$.
For any $\eta\in\Xi$, define
the corresponding vector-valued measure $\mu_\eta(\dd x) =
M^{-1}(\eta
)f(x) \eta(\dd x)$ and
note that $\mu_\eta\in\mathcal{M}_0$. We obtain
\begin{eqnarray*}
c^T D(\eta) c &=& c^T M^{-1} (\eta) B (\eta,
\eta) M^{-1} (\eta) c = \Phi_c (\mu_\eta)
\\
& \geq& \min_{\mu\in\mathcal{M}_0} \Phi_c (\mu)= \Phi_c (
\mu_0)= c^TD(\xi)c.
\end{eqnarray*}
Since the design $\xi$ does not depend on the particular vector $c$, it
follows that
$\xi$ is universally optimal.

Let us now prove (b) of Theorem~\ref{thuniversal-optimality}. Assume
that the design $\xi$ is universally optimal and let $g(x)$ be the function
associated with this design and computed by~(\ref{eqintKfxi=Lf+g}).

Consider first the case $m=1$.
In this case, the assumption that $\xi$ is universally optimal design
coincides with the assumption of simple optimality.
Also,
since $f(x) \neq0$ for all $x \in{\cal X}$, we can define $\gamma(x)
= g(x)/f(x) $ for all $x \in{\cal X}$.
In this notation, the statement (b) of Theorem \ref
{thuniversal-optimality} coincides with the statement of Corollary~\ref
{cor3.1}.

Assume now $m>1$.
Since the design $\xi$ is universally optimal it is $c$-optimal for
any vector $c$ and therefore
the necessary condition for $c$-optimality should be satisfied; this
condition is
\[
r_c(x,\xi)=c^T\bM^{ -1}(\xi)g(x)
f^T(x)\bM^{ -1}(\xi)c \geq0
\]
for all $x \in\mathcal{X}$ and $r_c(x,\xi)=0$ for all $x$ in the
support of the measure $\xi$. If
$g(x)=\gamma(x) f(x)$, where $\gamma(x)\geq0 $, then
\[
r_c(x,\xi)=\gamma(x) \bigl[ c^T\bM^{ -1}(
\xi)f(x) f^T(x)\bM^{
-1}(\xi)c \bigr] \geq0
\]
for any vector $c$ and all $x$ so that the necessary condition for
$c$-optimality is satisfied.
On the other hand, if $g(x)=\gamma(x) f(x)$, but $\gamma(x_0)< 0 $ for
some $x_0 \in\cal X$ then [in view of the fact that the matrix
$\bM^{ -1}(\xi)f(x_0) f^T(x_0)\bM^{ -1}(\xi)$ is nondegenerate] there
exists $c$ such that
$r_c(x_0,\xi)<0$ and the necessary condition for $c$-optimality of the
design $\xi$ is not satisfied.

Furthermore, if
the representation $g(x)=\gamma(x) f(x)$ does not hold, then
there exists a point $x_0 \in\mathcal{X}$ such that $g(x_0)\neq0$ and
$g(x_0)$ is not proportional to $f(x_0)$
[recall also that \mbox{$f(x)\neq0$} for all $x \in\mathcal{X}$].
Then
\[
r_c(x_0,\xi)=c^T\bM^{ -1}(
\xi)g(x_0)f^T(x_0)\bM^{ -1}(
\xi)c=c^Tab^Tc
\]
with $a=\bM^{ -1}(\xi)g(x_0)$ and $b=\bM^{ -1}(\xi)f(x_0)$.
Using Lemma~\ref{lemcabc<0},
we deduce that there exists a vector $c$ such that $r_c(x_0,\xi)<0$.
Therefore the design $\xi$ is not $c$-optimal and as a consequence also
not universally optimal.
\end{pf*}

In the one-parameter case ($m=1$) it is easy to construct examples
where the function $g(x)$ corresponding to the optimal design is
nonzero. For
example, consider the regression model $y(t)=\theta t + \varepsilon
(t)$, $t \in[-1,1]$, with the so-called spherical correlation function
\[
\rho(u)=1- \tfrac32 |u|/R + \tfrac12 \bigl(|u|/R\bigr)^3
\]
with $R=2$.
Then the design
assigning weights 0.5 to the points $-1$ and $1$ is optimal. For this
design, the function $g(x)$ defined in (\ref{eqintKfxi=Lf+g}) is
equal to
$g(x)= x(1-x^2)/16$, while the function $\gamma(x)$ is $\gamma
(x)=(1-x^2)/16$, $x \in[-1,1]$.

\section{Optimal designs for specific kernels and models}
\label{general}

\subsection{Optimality and Mercer's theorem}

In this section we consider the case
when the regression functions are proportional to eigenfunctions from
Mercer's theorem.
To be precise, let ${\cal X}$ denote a compact subset of a metric
space, and
let $\nu$ denote a measure on the corresponding Borel field with
positive density.
Consider the integral operator
%
%
\begin{equation}
\label{opdef} T_K(f) (\cdot) = \int_{\cal X} K(
\cdot,u)f(u) \nu(\dd u)
\end{equation}
on $L_2(\nu)$. Under certain assumptions on the kernel
[e.g., if $K(u,v)$ is symmetric, continuous and positive definite]
$T_K$ defines a symmetric,
compact self-adjoint operator.
In this case Mercer's theorem [see, e.g., \citet{kanwal1997}]
shows that there exist a countable number of eigenfunctions $\varphi_1,
\varphi_2, \ldots$ with positive eigenvalues
$\lambda_1,\lambda_2,\ldots$ of the operator $K$, that is,
%
%
\begin{equation}
\label{operator} T_k(\varphi_\ell) = \lambda_\ell
\varphi_\ell,\qquad \ell=1,2,\ldots.
\end{equation}
The next statement follows directly from Theorem~\ref{thuniversal-optimality}.

%
\begin{theorem}\label{thKL}
Let $\mathcal{X}$ be a compact subset of a metric space,
and assume that the covariance kernel
$K(x,u)$ defines an integral operator $T_K$ of the form~(\ref{opdef}),
where the eigenfunctions satisfy
(\ref{operator}).
Consider\vspace*{1pt} the regression model (\ref{lin-model}) with $f(x)=L (\varphi
_{i_1}(x),\ldots,\varphi_{i_m}(x))^T$
and the covariance kernel $K(x,u)$, where $L \in\mathbb\Real^{m\times
m}$ is a nonsingular matrix.
Then the design $\nu$ is universally optimal.
\end{theorem}

We note that the Mercer expansion is known analytically for certain
covariance kernels.
For example, if $\nu$ is the uniform distribution on the interval
$\mathcal{X}=[-1,1]$, and
the covariance kernel is of exponential type, that is,
$K(x,u)=e^{-\lambda|x-u|}$,
then the eigenfunctions are given by
\[
\varphi_k(x)=\sin(\omega_k x+k \pi/2),\qquad k \in\mathbb N,
\]
where $\omega_1,\omega_2,\ldots$ are positive roots of
the equation $\tan(2\omega)=-2\lambda\omega/(\lambda^2-\omega^2)$. Similarly,
consider as a second example, the covariance kernel $K(x,u)=\min\{x,u\}
$ and
$\mathcal{X}=[0,1]$, In this case, the eigenfunctions of the corresponding
integral operator are given by
\[
\varphi_k(x)=\sin\bigl((k+1/2)\pi x\bigr),\qquad k \in\mathbb N.
\]
In the following subsection we provide a further
example of the application of Mercer's theorem, which is of importance
for series estimation in nonparametric regression.

\subsection{Uniform design for periodic covariance functions}

Consider the regression functions
%
%
\begin{equation}
\label{eq3} f_j(x) = \cases{1, &\quad if $j=1$,
\cr
\sqrt{2} \cos
\bigl(2 \pi(j-1)x\bigr), &\quad if $j \geq2$,}
\end{equation}
and the design space $\mathcal{X}=[0,1]$. Linear models of the form
(\ref{lin-model}) with regression functions (\ref{eq3}) are widely
applied in series estimation of a nonparametric regression function
[see, e.g., Efromovich (\citeyear{efromovich1999,efromovich2008}) or
\citet{tsybakov2009}]. Assume that the correlation function
$\rho(x)$ is periodic with period $1$, that is, $\rho(x)=\rho(x+1)$,
and let a covariance kernel be defined by $ K(u,v) = \sigma^2\rho(u-v)
$ with $\sigma^2=1$. An example of the covariance kernel $\rho(x)$
satisfying this property is provided by a convex combination of the
functions $\{\cos(2\pi x), \cos^2(2\pi x),\ldots\}$.

%
\begin{theorem}
Consider regression model (\ref{lin-model}) with regression functions
$f_{i_1}(x), \ldots, f_{i_m}(x)$ $(1 \leq i_1 < \cdots< i_m)$ defined in
(\ref{eq3}) and
a correlation function $\rho(t)$ that is periodic with period $1$.
Then the uniform design is universally optimal.
\end{theorem}
\begin{pf}
We will show that the identity
%
%
\begin{equation}
\label{eq4} \int^1_0 K(u,v) f_j(u)
\,\dd u=\int^1_0 \rho(u-v)f_j(u)\,
\dd u= \lambda_j f_j(v)
\end{equation}
holds for all $v\in[0,1]$,
where $\lambda_j=\int\rho(u)f_j(u)\,\dd u$ ($j \ge1$). The assertion
then follows from Theorem~\ref{thKL}.

To prove (\ref{eq4}), we define
$
A_j(v) = \int^1_0 \rho(u-v) f_j(u) \,\dd u
$
which should be shown to be $\lambda_j f_j(v) $.
For $j=1$ we have $A_1(v)=\lambda_1$ because
$
\int^1_0 \rho(u-v)\,\dd u= \int^1_0 \rho(u)\,\dd u =\lambda_1
$
by the periodicity of the function $\rho(x)$. For $j=2,3,\ldots$ we
note that
\begin{eqnarray*}
A_j(v) &=& \int^1_0
\rho(u-v)f_j(u)\,\dd u = \int^{1-v}_{-v}
f_j(u+v)\rho(u)\,\dd u
\\
&=& \int^{1-v}_0 f_j(u+v)\rho(u)\,
\dd u + \int^0_{-v}f_j(u+v)\rho(u)\,
\dd u.
\end{eqnarray*}
Because of the periodicity we have
\[
\int^0_{-v} f_j (u+v) \rho(u)\,\dd
u = \int^1_{1-v} f_j(u+v)\rho(u)\,
\dd u,
\]
which gives
$
A_j(v)= \int^1_0 f_j (u+v) \rho(u) \,\dd u.
$
A simple calculation now shows
%
%
\begin{equation}
\label{eq6} A^{\prime\prime}_j (v)=-b^2_j
A_j(v),
\end{equation}
where $b^2_j = (2 \pi(j-1))^2$ and
\begin{eqnarray*}
A_j (0) &=& \sqrt2 \int^1_0 \cos
\bigl(2 \pi(j-1)u\bigr) \rho(u)\,\dd u = \sqrt{2} \lambda_j,
\\
A^{\prime}_j (0) &=& -b_j \sqrt2 \int
^1_0 \sin\bigl(2 \pi(j-1)u\bigr) \rho(u)\,\dd u
= 0.
\end{eqnarray*}
Therefore (from the theory of differential equations) the unique
solution of (\ref{eq6}) is of the form
$
A_j(v)=c_1 \cos(b_jv)+c_2 \sin(b_jv),
$
where $c_1$ and $c_2$ are determined by initial conditions, that is,
$
A_j(0) = c_1 = \sqrt{2} \lambda_j,
A^\prime_j (0) = b_j c_2 = 0.
$
This yields
$
A_j(v) = \lambda_j \sqrt{2} \cos(2 \pi(j-1)v) = \lambda_j f_j(v)
$
and proves identity~(\ref{eq4}).
\end{pf}

%

\subsection{Optimal designs for the triangular covariance function}
\label{triang}

Let us now consider the triangular correlation function defined by
%
%
\begin{equation}
\label{eqtriang} \rho(x)=\max\bigl\{0,1-\lambda|x|\bigr\}.\vadjust{\goodbreak}
\end{equation}
On the one hand this function arises as a correlation function of the
process of increments of a Brownian motion
which is in turn related to a Brownian bridge
after a suitable conditioning is made; see \citet{mehmcf1965}.
On the
other hand it is motivated by the fact that for ``small'' values of the
parameter
$\lambda$, it provides a good approximation of the exponential
correlation kernel $\rho_\lambda(x)=\exp(-\lambda|x|)$, which
is widely used for modeling correlations in regression models; see
\citet
{uciatk2004} or \citet{detpephol2010}, among others. For the
exponential correlation kernel optimal designs are difficult to find,
even in the linear regression
model; see \citet{detkunpep2008}. However, as the next theorem
shows, it
is possible to explicitly derive optimal designs for the linear model
with a triangular correlation function. It will be demonstrated in
Example~\ref{exam3}
below that for ``small'' and ``moderate'' values of the parameter
$\lambda$, these designs provide also an efficient solution of the
design problem for the exponential
correlation kernel.

%
\begin{theorem}\label{THtriang}
Consider model (\ref{lin-model}) with $f(x)=(1,x)^T$, $\mathcal
{X}=[-1,1]$ and
the triangular correlation function (\ref{eqtriang}).

\begin{longlist}[(a)]
\item[(a)]
If $\lambda\in(0,1/2]$, then the design $\xi^*=\{-1,1;1/2,1/2\}$
is universally optimal.

\item[(b)] If $\lambda\in\mathbb{N}$, then
the design supported at $2\lambda+1$ points
$x_k=-1+k/\lambda$, $k=0,1,\ldots,2\lambda$, with equal weights is
universally optimal.
\end{longlist}
\end{theorem}

\begin{pf}
To prove part (a) we will show that
$\int\rho(x-u)f_i(u)\xi^*(\dd u)=\lambda_i f_i(x)$ for $i=1,2$ and some
$\lambda_1$ and $\lambda_2$.
By direct calculations we obtain for $f_1(x)=1$,
\[
\int^1_{-1}\rho(x-u)\xi^*(\dd u)=
\frac{-\lambda|x + 1|}{2}+\frac
{1 -
\lambda|x - 1|}{2}=1-\lambda
\]
and, consequently, $\lambda_1=1-\lambda$.
Similarly, we have for $f_2(x)=x$,
\[
\int^1_{-1} u\rho(x-u)\xi^*(\dd u) = -
\frac{1 - \lambda|x +
1|}2+\frac
{1 - \lambda|x - 1|}2=\lambda x
\]
and, therefore, $\lambda_2=\lambda$.
Thus, the assumptions of Theorem~\ref{thuniversal-optimality} are fulfilled.

Part (b). Straightforward but tedious calculations show that
$M(\xi^*)=\operatorname{diag}(1,\gamma)$,
where $\gamma=\sum_{k=0}^{2\lambda+1} x_k^2/(2\lambda+1)=(\lambda
+1)/(2\lambda)$.
Also we have
\[
\int\rho(x-u)f_i(u)\xi^*(\dd u)=\lambda_i f_i(x)
\]
for $i=1,2$
where $\lambda_1=\lambda_2=1/(2\lambda+1)$.
Thus, the assumptions of Theorem~\ref{thuniversal-optimality} are fulfilled.
\end{pf}

The designs provided in Theorem~\ref{THtriang} are also optimal for
the location scale model; see \citet{zhidetpep2010}.
However, unlike the results of previous subsections
the result of Theorem~\ref{THtriang} cannot be extended to polynomial
models of higher order.

We conclude this section with an example illustrating the efficiency of
the designs for the triangular kernel in models
with correlation structure defined by the exponential kernel.

%
\begin{example} \label{exam3}
Consider the location scale model [$f(x)=1$] and the linear regression
model [$f(x)=(1,x)^T$], $\mathcal{X}=[-1,1]$, and the correlation
function $\rho(x)=\exp\{- \lambda|x|\}$. In Table~\ref{tabtrin} we
display the $D$-efficiencies of the universally optimal design
%
%
\begin{table}
\tablewidth=260pt
\caption{$D$-Efficiencies of the universally optimal design $\xi=\{
-1,1;\break0.5,0.5\}$
calculated under the assumption of a triangular correlation function in
the constant and linear
regression model with the exponential correlation function $\rho
(x)=e^{-\lambda|x|}$}
\label{tabtrin}
%
\begin{tabular*}{\tablewidth}{@{\extracolsep{\fill}}lc ccc c@{}}
\hline
$\bolds{\lambda}$&\textbf{0.1}& \textbf{0.3}& \textbf{0.5}& \textbf{0.7}& \textbf{0.9}\\
\hline
Constant &0.999&0.997&0.978&0.946&0.905\\
Linear & 0.999&0.999&0.991&0.974&0.950\\
\hline
\end{tabular*}
\end{table}
calculated under the assumption of the triangular kernel
(\ref{eqtriang}) for various values of the parameter
$\lambda\in[0.1,0.9]$. For this design we observe in all cases a
$D$-efficiency of at least $90\%$. In most cases it is higher than
$95\%$.
\end{example}

\subsection{Polynomial regression models and singular kernels}
\label{singular}

In this section we consider the polynomial regression model, that is,
$f(x) = (1,x,\ldots, x^{m-1})^T$,
with logarithmic covariance kernel
%
%
\begin{equation}
\label{logar} K(u,v) = \gamma-\beta\ln(u-v)^2,\qquad \beta>0, \gamma\geq0,
\end{equation}
and the kernel
%
%
\begin{equation}
\label{betakern} K(u,v)=\gamma+\beta/|u-v|^\alpha,\qquad 0 \leq\alpha< 1,
\gamma\geq0, \beta>0,
\end{equation}
for which the universally optimal designs can be found explicitly.

Covariance functions ${K}(u,v)$ with a singularity at $u=v$ appear
naturally as approximations to many
standard covariance functions $\tilde{K}(u,v)=\sigma^2 \tilde\rho
(u-v)$ with $\tilde\rho(0)=1$ if
$\sigma^2$ is large. A general scheme for this type of approximation is
investigated in \citet{zhidetpep2010}, Section 4. More precisely, these
authors discussed the case where the covariance kernel can be
represented as
$\sigma^2_\delta\tilde\rho_\delta(t)= r \ast h_{\delta} (t)$
with a singular kernel $r(t)$ and a smoothing kernel $h_{\delta}(\cdot
)$ (here $\delta$ is a smoothing parameter and $\ast$ denotes the
convolution operator). The basic idea is illustrated in the following example.
%
%
\begin{example} \label{exam4}
Consider the covariance kernel $\tilde{K}(u,v)= \rho_\delta(u-v)$,
where
%
%
\begin{equation}
\label{covcon} \rho_\delta(t) = 2 - {1\over\delta} \log\biggl(
{ |t+\delta|^{t+\delta} \over
|t-\delta
|^{t-\delta} } \biggr).
\end{equation}
For several values of $\delta$, the function $\rho_\delta$ is displayed
in Figure~\ref{figsm-k}.
A straightforward calculation shows that
$
\rho_\delta(t) =r \ast h_{\delta} (t),
$
where $r(t)=-\ln(t)^2$ and $h_{\delta}$ is the density of the uniform
distribution on the interval
on $[-\delta,\delta]$. As illustrated by Figure~\ref{figsm-k}, the
function $\rho_\delta(\cdot)$ is well approximated
by the singular kernel $r(\cdot)$ if $\delta$ is small.

%
%
\begin{figure}

\includegraphics{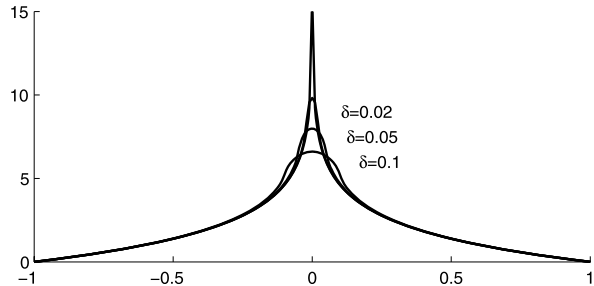}

\caption{The logarithmic covariance kernel $r(t)=-\ln(t)^2$ and
the covariance kernel (\protect\ref{covcon}), where $\delta=0.02,0.05,0.1$.}
\label{figsm-k}
\end{figure}

In\vspace*{1pt} Figure~\ref{figdopt-log} we display the $D$-optimal designs
(constructed numerically) for the quadratic model with
a stationary error process with covariance kernel
$ \tilde{K}(u,u+t)=\rho_\delta(t)$, where $\rho_\delta$ is defined in
(\ref{covcon}) and
$\delta=0.02,0.05,0.1$.
As one can see, for small $\delta$ these designs are very close
to the arcsine design, which is the $D$-optimal design for the
quadratic model and
the logarithmic kernel, as proved in Theorem~\ref
{THmulti-genarcsine45} of the following section.

%
%
\begin{figure}[b]

\includegraphics{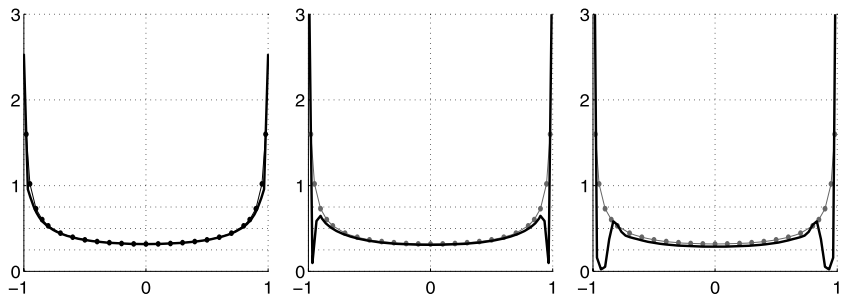}

\caption{Density functions corresponding to the $D$-optimal designs
for the quadratic model with covariance kernel (\protect\ref
{covcon}), where
$\delta=0.02$ (left), $\delta=0.05$ (middle) and $\delta=0.1$ (right).
The $y$-axis corresponds\vspace*{1pt} to values of the density functions. The
corresponding designs are obtained by (\protect\ref{despoint}), where
$a^{-1}$ is the distribution function corresponding to the displayed
densities. The grey\vspace*{1pt} line corresponds to the arcsine density
$p(x)={1}/{\pi\sqrt{1-x^2}}$.} \label{figdopt-log}
\end{figure}

In Table~\ref{tabeff-logd} we show the efficiency of the arcsine
distribution [obtained by maximizing
$\det(D(\xi) )$ with the logarithmic kernel] in the quadratic
regression model with the kernel
(\ref{covcon}). We observe a very high efficiency with respect to the
$D$-optimality criterion. Even in the case $\delta=0.1$ the
efficiency is $93.6\%$ and it converges quickly to $100\%$ as $\delta$
approaches $0$.
\end{example}
%
%
\begin{example} \label{exam4a}
The arcsine density can also be used as an alternative approximation
to the
exponential correlation function or correlation functions of a similar
type, that is,\vspace*{1pt}
$ \rho_{\lambda,\nu} (t)=\exp({-\lambda|t|^{\nu} } ) $.
For the case $\lambda=\nu=1$ the function
$\frac7{15} (1-\frac3{17} \ln t^2 )$ can be considered as a
reasonable approximation to $ \exp( {-|t|})$
on the interval $[-1,1]$; see the left part of Figure
\ref{figexpacsin}. Similarly, if $\lambda=1$, $\nu=1/4$, it is
illustrated in the right part of Figure
\ref{figexpacsin} that the function $\frac38 -\frac1 {25} \ln t^2 $
provides a very accurate approximation
of the exponential correlation function.
%
%
\begin{table}
\tablewidth=250pt
\caption{Efficiency of the arcsine design $\xi_{a}$
for the quadratic model~and~the kernel (\protect\ref{covcon})}
\label{tabeff-logd}
\begin{tabular*}{\tablewidth}{@{\extracolsep{\fill}}l ccc cc@{}}
\hline
$\delta$&0.02& 0.04& 0.06& 0.08& 0.1\\
[4pt]
Eff$(\xi_a)$ &0.998&0.978&0.966&0.949&0.936\\
\hline
\end{tabular*}
\end{table}
%
\begin{figure}[b]

\includegraphics{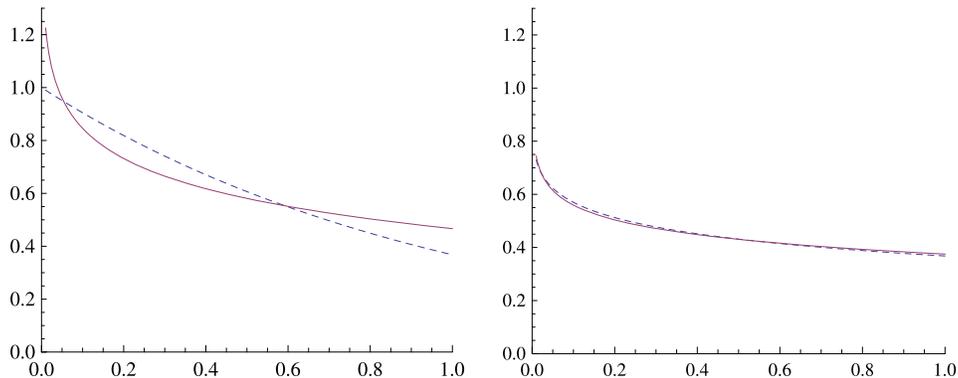}

\caption{Left panel: the function $\frac7{15}(1-\frac3{17} \ln t^2 )$
(solid line) as
an approximation of the exponential correlation function $\exp(
{-|t|})$. Right panel: the function
$\frac38 -\frac1 {25} \ln t^2 $ (solid line) as
an approximation of the exponential correlation function $\exp(
-|t|^{1/4})$ (dashed line).}
\label{figexpacsin}
\end{figure}
As a consequence, the arcsine design (optimal for the logarithmic
kernel) will also have a high efficiency with respect to
these kernels, and this argument is illustrated in Table \ref
{tabeff-exp} of Section~\ref{num1} where we calculate
the $D$-efficiencies of the arcsine design in polynomial regression
models with correlation function $ \exp( {-|t|})$.
For the correlation function $ \exp({- \lambda|t|^{1/4} )}$ a similar
$D$-efficiency of the arcsine design can be observed.
For example, if $\lambda=0.5$, $2.5$ the $D$-efficiencies of the
arcsine design in the linear regression model
are $100\%$ and $96.9\%$, respectively, while they are $99.9\%$ and
$97.1\%$ in the quadratic model. Other
choices of $\lambda$ and $\nu$ yield similar results, which are not
displayed for the sake of brevity.
\end{example}

\subsubsection{Optimality of the arcsine design}

We will need the following lemma, which states a result in the theory
of Fredholm--Volterra integral equations; see \citet{mashan2002},
Chapter 9, page 211.
%
%
\begin{lemma}\label{LemFredVoltEq-Cheb}
The Chebyshev polynomials of the first kind
$T_n(x)=\cos(n\* \arccos x)$ are the eigenfunctions of the integral
operator with the kernel\break
$H(x,v)= - \ln(x-v)^2 / \sqrt{1-v^2}$. More precisely, for all $n=0,1,
\ldots$ we have for all $n \in\mathbb{N}$
\[
\lambda_n T_n(x) =-\int^1_{-1}
T_n(v) \ln(x-v)^2 \,\frac{\dd v}{\pi
\sqrt{1-v^2}},\qquad x\in[-1,1],
\]
where $\lambda_0=2\ln2$ and $\lambda_n=2/n$ for $n\geq1$.
\end{lemma}
With the next result we address the problem of uniqueness of the
optimal design. In particular, we
give a new characterization of the arcsine distribution. A proof can
be found in the \hyperref[app]{Appendix}.

%
\begin{theorem}\label{tharcsine}
Let $n$ be a nonnegative integer and $\zeta$ be a random variable
supported on the interval $[-1,1]$.
Then the distribution of $\zeta$ has the arcsine density (\ref{arcs})
if and only if the equality
\[
\mathbb{E} T_n(\zeta) \bigl(-\ln(\zeta-x)^2
\bigr)=c_nT_n(x)
\]
holds for almost all $x\in[-1,1]$, where $c_n=2/n$ if $n\in\mathbb{N}$
and $c_0=2\ln2$ if $n=0$.
\end{theorem}

The following result is an immediate consequence of Theorems \ref
{thuniversal-optimality} and~\ref{tharcsine}.

%
\begin{theorem}\label{THmulti-genarcsine45}
Consider the polynomial regression model (\ref{lin-model})
with $f(x)=(1, x, x^2,\ldots,x^{m-1})^T$, $x \in[-1,1]$,
and the covariance kernel (\ref{logar}). Then the arcsine design $\xi
_a$ with density (\ref{arcs})
is the universally optimal design.
\end{theorem}
\begin{pf}
We assume without loss of generality that $\beta=1$ and consider the
function $\rho(x)=-\ln x^2+\gamma$ with positive $\gamma$. From Lemma
\ref{LemFredVoltEq-Cheb} we obtain
\begin{eqnarray*}
\int^1_{-1} \bigl(-\ln(u-x)^2+\gamma
\bigr) T_n(u)p(u)\,\dd u &=& - \int^1_{-1}
\ln(u-x)^2 T_n(u)p(u)\,\dd u
\\
&=& \lambda_n T_n(x) + \gamma\delta_{n0},
\end{eqnarray*}
where $\delta_{xy}$ denotes Kronecker's symbol and we have used the
fact that
$\int_{-1}^1 T_n(u)/ \sqrt{1 - u^2} \,\dd u=0$ whenever $n \geq1$.
This proves\vadjust{\goodbreak} (\ref{eqintKfxi=Lf+g}) for the arcsine distribution
and
the vector $t(x)=(T_0(x),\ldots, T_{m-1}(x))^T$ where the function
$g(x)$ is equal to~0 for all $x$.
Now $f(x)=(1,x,\ldots,x^{m-1})^T = Lt(x)$ for some nonsingular $m
\times m$ matrix. Therefore (\ref{eqintKfxi=Lf+g}) holds also for the
vector $f(x)$ with $g(x) \equiv0$ (and a different matrix~$\Lambda$).
The statement of the theorem now follows from Theorems \ref
{thuniversal-optimality} and~\ref{tharcsine}.
\end{pf}

\subsubsection{Generalized arcsine designs}

For $\alpha\in(0,1)$
consider the Gegenbauer polynomials $C^{(\alpha)}_m(x)$
which are orthogonal with respect to the weight function
%
%
\begin{equation}
\label{weightg} p_\alpha(x)= \frac{(\Gamma(\alpha+
{1/2}))^2}{2^\alpha
\Gamma
(2 \alpha+1 )} \bigl(1-x^2
\bigr)^{(\alpha-1)/2},\qquad x \in[-1,1].
\end{equation}
For the choice $\alpha=0$ the Gegenbauer polynomials $C^{(\alpha)}_m(x)$
are proportional to the Chebyshev polynomials of the first kind
$T_m(x)$. Throughout this paper we will call
the corresponding beta-distributions \textit{generalized arcsine designs}
emphasizing the fact that
the distribution is symmetric and the parameter $\alpha$ varies in the
interval $(0,1)$.
The following result [from the theory of Fredholm--Volterra integral equations
of the first kind with special kernel,
see \citet{fahabddar1999}] establishes an analogue of
Lemma~\ref{LemFredVoltEq-Cheb} for the kernel
%
%
\begin{equation}
\label{gegkern} H(u,v)= \frac{1}{|u-v|^\alpha(1-v^2)^{(1- \alpha)/2}}.
\end{equation}

%
\begin{lemma}\label{LemFredVoltEq-Geg}
The Gegenbauer polynomials $C^{(\alpha/2)}_n(x)$ are the
eigenfunctions of the integral operator with the kernel defined in
(\ref{gegkern}).
More precisely, for all $n=0,1, \ldots$ we have
\[
\lambda_n C^{(\alpha/2)}_n(x) =-\int
^1_{-1} \,\frac{1}{|x-v|^\alpha} C^{(\alpha/2)}_n(v)
\frac{\dd v}{ (1-v^2)^{(1-\alpha)/2} }
\]
for all $x\in[-1,1]$, where $\lambda_n=\frac{\pi\Gamma(n+\alpha
)}{\cos(\alpha\pi/2)\Gamma(\alpha)n!}$.
\end{lemma}

The following result generalizes Theorem 8 of \citet{zhidetpep2010}
from the case of a location scale model to polynomial regression models.


%
\begin{theorem}\label{THmulti-genarcsine}
Consider the polynomial regression model (\ref{lin-model})
with $f(x)=(1, x, x^2,\ldots,x^{m-1})^T$, $x \in[-1,1]$,
and covariance kernel (\ref{betakern}).
Then the design with generalized arcsine density defined in (\ref{weightg})
is universally optimal.
\end{theorem}
\begin{pf}
It is easy to see that the optimal design does not depend on $\beta$
and we thus assume that $\beta=1$ in (\ref{betakern}).\vadjust{\goodbreak}

To prove the statement for the kernel $\rho(x)=1/|x|^\alpha+\gamma$
with positive $\gamma$
we recall the definition of $p_\alpha$ in (\ref{weightg}) and obtain
from Lemma~\ref{LemFredVoltEq-Geg}
\begin{eqnarray*}
\int\biggl(\frac{1}{|u-x|^\alpha}+\gamma\biggr) C^{({\alpha
/2})}_n(u)p_{{\alpha}/{2}}(u)
\,\dd u&=& \int\frac{1}{|u-x|^\alpha} C^{({\alpha/2})}_n(u)p_{
{\alpha
/2}}(u)
\,\dd u\\
&\propto& C^{({\alpha/2})}_n(x) 
\end{eqnarray*}
for any $n \in\mathbb{N}$ since $\int C^{(\alpha/2)}_n(u)p_{\alpha
/2}(u)\,\dd u=0$.
Consider the design $\xi$ with density~$p_{{\alpha}/{2}}$. For this
design, the function $g(x)$ defined in (\ref{eqintKfxi=Lf+g}) is
identically zero; this follows from the formula above.
It now follows by the same arguments as given at the end of the proof
of Theorem~\ref{THmulti-genarcsine45} that the design with density
$p_{{\alpha}/{2}}$ is universally optimal.
\end{pf}

\section{Numerical construction of optimal designs}
\label{numerics}

\subsection{An algorithm for computing optimal designs}

Numerical computation of optimal designs for a common linear regression
model (\ref{lin-model})
with given correlation function
can be performed by an extension of the multiplicative algorithm
proposed by \citet{detpepzhi2008} for the case of noncorrelated
observations.
Note that the proposed algorithm constructs a discrete design which can
be considered as an approximation to a design which satisfies the necessary
conditions of optimality of Theorem~\ref{thphi-opt-cond-pdim}. By
choosing a fine discretization $\{x_1,\ldots,x_n\}$
of the design space $\mathcal{X}$ and running the algorithm long enough,
the accuracy of approximation can be made arbitrarily small (in the
case when convergence is achieved).

Denote by $\xi^{(r)}=\{x_1,\ldots,x_n;w^{(r)}_1,\ldots,w^{(r)}_n\}$ the
design at the iteration $r$,
where
$w^{(0)}_1,\ldots,w^{(0)}_n$ are nonzero weights, for example, uniform.
We propose the following updating rule for the weights:
%
%
\begin{equation}
\label{equpd-rule} w^{(r+1)}_i=\frac{w^{(r)}_i (\psi(x_i,\xi
^{(r)})-\beta_r
)}{\sum_{j=1}^n w^{(r)}_j (\psi(x_j,\xi^{(r)})-\beta_r )},\qquad i=1,
\ldots,n,
\end{equation}
where $\beta_r$ is a tuning parameter [the only condition on $\beta_r$
is the positivity of all the weights in (\ref{equpd-rule})],
$
\psi(x,\xi)={\varphi(x,\xi)}/{b(x,\xi)}
$
and the functions $\varphi(x,\xi)$ and $b(x,\xi)$ are defined in
(\ref{eqphi-x-xi})
and (\ref{eqb-x-xi}), respectively.
Condition (\ref{eqphi-ekv-multi}) takes the form $
\psi(x,\xi^*)\le1
$ for all $x \in\mathcal{X}$.
Rule (\ref{equpd-rule}) means that at the next iteration the weight of
a point $x=x_j$ increases if condition (\ref{eqphi-ekv-multi}) does
not hold at this point.

A measure $\xi_*$ is a fixed point of the iteration (\ref{equpd-rule})
if and only if
\mbox{$
\psi(x,\xi_*) = 1
$}
for all $x \in\operatorname{supp} (\xi_*)$ and $
\psi(x,\xi_*) \leq1
$
for all $x \in\mathcal{X} \setminus\operatorname{supp} (\xi_*)$. That
is, a~design $\xi_*$ is a
fixed point of the iteration (\ref{equpd-rule}) if and only if it
satisfies the optimality condition of Theorem~\ref{thphi-opt-cond-pdim}.
We were not able to theoretically prove the convergence of iterations
(\ref{equpd-rule}) to the design satisfying
the optimality condition of Theorem~\ref{thphi-opt-cond-pdim}, but we
observed this convergence in all numerical studies.
In particular, for the cases where we could derive the optimal designs
explicitly, we observed convergence of the algorithm to the optimal design.


Algorithm (\ref{equpd-rule}) can be easily extended to cover the case
of singular covariance kernels.
Alternatively, a singular kernel can be approximated by a nonsingular
one using the technique
described in \citet{zhidetpep2010}, Section 4.


\subsection{Efficiencies of the uniform and arcsine densities} \label{num1}

In the present section we numerically study the efficiency (with
respect to the $D$-optimality criterion) of the uniform and arcsine
designs for the polynomial model (\ref{lin-model}) with
$f(x)=(1,x,\ldots,x^{m-1})^T$ and the exponential correlation function
$\rho(t)=e^{-\lambda|t|}$, $t \in[-1,1]$. We determine the efficiency
of a design $\xi$ as
\[
\operatorname{Eff}(\xi)= \biggl(\frac{\det\bD(\xi^*)}{\det\bD(\xi
)} \biggr)^{1/m},
\]
where $\xi^*$ is the design computed by the algorithm described in
the previous section (applied to the $D$-optimality criterion). The
results are depicted in Table~\ref{tabeff-exp}.
%
%
\begin{table}[b]
\tablewidth=278pt
\caption{Efficiencies of the uniform design $\xi_u$ and the arcsine
design $\xi_{a}$
for~the~polynomial regression model of degree $m-1$ and
the~exponential~correlation function $\rho(x)=e^{-\lambda|x|}$}
\label{tabeff-exp}
%
\begin{tabular*}{\tablewidth}{@{\extracolsep{\fill}}lc ccc ccc@{}}
\hline
&$\bolds{\lambda}$&\textbf{0.5}& \textbf{1.5}& \textbf{2.5}& \textbf{3.5}& \textbf{4.5}
& \textbf{5.5}\\
\hline
$m=1$&Eff$(\xi_u)$&0.913&0.888&0.903&0.919&0.933&0.944\\
&Eff$(\xi_a)$&0.966&0.979&0.987&0.980&0.968&0.954\\
[4pt]
$m=2$&Eff$(\xi_u)$&0.857&0.832&0.847&0.867&0.886&0.901\\
&Eff$(\xi_a)$&0.942&0.954&0.970&0.975&0.973&0.966\\
[4pt]
$m=3$&Eff$(\xi_u)$&0.832&0.816&0.826&0.842&0.860&0.876\\
&Eff$(\xi_a)$&0.934&0.938&0.954&0.968&0.976&0.981\\
[4pt]
$m=4$&Eff$(\xi_u)$&0.826&0.818&0.823&0.835&0.849&0.864\\
&Eff$(\xi_a)$&0.934&0.936&0.945&0.957&0.967&0.975\\
\hline
\end{tabular*}
\end{table}
We observe that
the efficiency of the arcsine design is always higher than the
efficiency of the uniform design.
Moreover, the absolute difference between the efficiencies of the two designs
increases as $m$ increases.
On the other hand, in most cases the efficiency of the uniform design
and the arcsine design decreases as $m$ increases.

\section{Conclusions}
\label{concl}

In this paper we have addressed the problem of constructing optimal
designs for
least squares estimation in regression models with correlated observations.
The main challenge in problems of this\vadjust{\goodbreak} type is that---in contrast to
``classical'' optimal design theory for uncorrelated data---the
corresponding optimality
criteria are not convex (except for the location scale model).
By relating the design problem to an integral operator problem,
universally optimal design can be identified explicitly for a broad
class of regression
models and correlation structures. Particular attention is paid
to a trigonometric regression model involving only cosines terms,
where it is proved that the uniform distribution is universally optimal
for any periodic
kernel of the form $K(u,v)=\rho(u-v)$. For the classical polynomial
regression model with
a covariance kernel given by the logarithmic potential
it is proved that the arcsine distribution is universally optimal.
Moreover, optimal designs are derived for several other regression models.

So far
optimal designs for regression models with correlated observations have only
be derived explicitly for the location scale model, and to our best knowledge
the results presented in this paper provide the first explicit
solutions to this
type of problem for a general class of models with more than one parameter.

We have concentrated on the construction of optimal
designs for least squares estimation (LSE) because the best linear unbiased
estimator (BLUE) requires the knowledge of the correlation matrix.
While the BLUE is often sensitive with respect to misspecification of the
correlation structure, the corresponding optimal designs for the LSE
show a remarkable
robustness. 
Moreover, the difference between BLUE and LSE is often surprisingly
small, and in
many cases BLUE and LSE with certain correlation functions
are asymptotically equivalent; see \citet{rao1967}, \citet
{kruskal1968}.

Indeed, consider the location scale model $ y(x)=\theta+\ve(x) $ with
$K(u,v)=\rho(u-v)$, where the knowledge of a full trajectory of a
process $y(x)$ is available. Define the (linear unbiased) estimate $
\hat\theta(G)=\int y(x) \,\dd G(x), $ where $G(x)$ is a distribution
function of a signed probability measure. A celebrated result of
\citet{grenander1950} states that the ``estimator''
$\hat\theta(G^*)$ is BLUE if and only if $\int\rho(u-x)\,\dd G^*(u)$ is
constant for all $x\in\mathcal{X}$. This result was extended by
N{\"a}ther [(\citeyear{naether1985a}), Section~4.3], to the case of
random fields with constant mean. Consequently, if $G^*(x)$ is a
distribution function of a nonsigned (rather than signed) probability
measure, then LSE coincides with BLUE and an asymptotic optimal design
for LSE is also an asymptotic optimal design for BLUE.
\citet{hajek1956} proved that $G^*$ is a distribution function of
a nonsigned probability measure if the correlation function $\rho$ is
convex on the interval $(0,\infty)$. \citet{zhidetpep2010} showed
that $G^*$ is a proper distribution function for a certain families of
correlation functions including nonconvex ones.

In Theorem~\ref{thuniversal-optimality} we have characterized the cases
where there exist universally optimal designs for
ordinary least squares estimation.
Specifically, a design $\xi^*$ is universally optimal for least squares
estimation
if and only if condition (\ref{eqintKfxi=Lf+g}) with $g(x)\equiv0$ is
satisfied.
Moreover, the proof of Theorem~\ref{thuniversal-optimality} shows that
in this case
the signed vector-valued measure
\[
\mu(\dd x)=M^{-1}\bigl(\xi^*\bigr)f(x)\xi^*(\dd x)
\]
and the LSE
minimizes (with respect to the Loewner ordering)
the matrix
\[
\iint K(x,u)\mu(\dd x)\mu^T(\dd u)
\]
in the space $\mathcal{M}$
of all vector-valued signed measures. Because this matrix is the
covariance of the linear estimate $\int y(x) \mu(\dd x)$ (where $\mu$ is
a vector of signed measures) it follows
that under the assumptions of Theorem~\ref{thuniversal-optimality},
the LSE combined with the universally optimal design $\xi^*$ give
exactly the same asymptotic covariance matrix as the BLUE and the
optimal design for the BLUE.

%
\begin{appendix}\label{app}
\section*{Appendix: Some technical details} 

\begin{pf*}{Proof of Lemma~\ref{lem-convex}}
For any $c\in\mathbb{R}^m$ and $\mu\in\mathcal{M}$ we set $\nu
(\cdot
) = c^T \mu(\cdot)$, where $\nu(\dd x)$ is a signed measure on
$\mathcal{X}$.
Then the functional
\[
\Phi_c(\mu) = c^T \iint K(x,u) \mu(\dd x)
\mu^T (\dd u) c
\]
can also be written as
\[
\Phi_c(\mu)= \Psi(\nu)= \iint K(x,u) \nu(\dd x) \nu(\dd u).
\]

For any $\alpha\in[0,1]$ and any two signed measures $\nu_0$ and
$\nu_1$ on $\mathcal{X}$ we have
\begin{eqnarray*}
& &\Psi\bigl(\alpha\nu_0+(1-\alpha)\nu_1\bigr)
\\
&&\qquad= \iint K(u,v)\bigl[\alpha\nu_0(\dd u)+(1-\alpha)
\nu_1(\dd u)\bigr] \bigl[\alpha\nu_0(\dd v)+(1-\alpha)
\nu_1(\dd v)\bigr]
\\
&&\qquad= \alpha^2\iint K(u,v)\nu_0(\dd u)
\nu_0(\dd v)+(1-\alpha)^2\iint K(u,v)
\nu_1(\dd u)\nu_1(\dd v)
\\
&&\qquad\quad{}+ 2\alpha(1-\alpha)\iint K(u,v)\nu_0(\dd u)
\nu_1(\dd v)
\\
&&\qquad= \alpha^2\Psi(\nu_0)+(1-\alpha)^2
\Psi(\nu_1)+2\alpha(1-\alpha)\iint K(u,v)\nu_0(\dd u)
\nu_1(\dd v)
\\
&&\qquad= \alpha\Psi(\nu_0)+(1-\alpha)\Psi(\nu_1)-
\alpha(1-\alpha)A,
\end{eqnarray*}
where
\begin{eqnarray*}
A&=&\iint K(u,v)\bigl[\nu_0(\dd u)\nu_0(\dd v)+
\nu_1(\dd u)\nu_1(\dd v)-2\nu_0(\dd u)
\nu_1(\dd v)\bigr]
\\
&=&\iint K(u,v)\zeta(\dd u)\zeta(\dd v)
\end{eqnarray*}
and
$\zeta(\dd u)=\nu_0(\dd u)-\nu_1(\dd u)$.
In view of (\ref{eqposdef}),
we have $A\ge0$ and therefore the functional $\Psi(\cdot)$ is convex.
\end{pf*}

\begin{pf*}{Proof of Lemma~\ref{lemcabc<0}}
As vectors $a$ and $b$ are linearly independent, we have
$a^T a >0$, $b^T b >0$ and $(a')^T b' <1$, where $a'=a/\sqrt{a^T a }$
and $b'=b/\sqrt{ b^T b }$.
For any vector $c \in\mathbb{R}^m$, we can represent $S_c$ as
\[
S_c=c^Ta b^Tc = \sqrt{a^T a
\cdot b^T b } \cdot c^Ta' \cdot
c^T b'.
\]
With the choice $c=a'-b'$ it follows
\[
c^Ta'= 1- \bigl(a'\bigr)^T
b' >0 \quad\mbox{and}\quad c^Tb'=
\bigl(a'\bigr)^T b'-1 <0
\]
implying $S_c<0$.
\end{pf*}

%
%


\begin{pf*}{Proof of Theorem~\ref{tharcsine}}
Note that the part ``if'' of the statement follows from Lemma~\ref
{LemFredVoltEq-Cheb},
and we should prove the part ``only if.'' Nevertheless, we provide a
proof of the part ``if''
since it will be the base for proving the part ``only if.''

Since the statement for $n=0$ is proved in \citet{schzhi2009}, we
consider the case $n\in\mathbb{N}$ in the rest of proof.
Using the transformation
$\varphi=\arccos u$ and $\psi=\arccos x$, we obtain $T_n(\cos\varphi
)=\cos(n \varphi)$ and
\[
\int_{-1}^1 \frac{\ln(u-x)^2}{\pi\sqrt{1-u^2}}T_n(u)
\,\dd u= \int_0^\pi\frac{\ln(\cos\varphi-x)^2}{\pi\sin\varphi}\cos(n
\varphi)\sin\varphi\,\dd\varphi.
\]
Consequently, in order to prove Theorem~\ref{tharcsine} we have to
show that the function
\[
\int_0^\pi{\ln(\cos\varphi-\cos
\psi)^2}\cos(n\varphi) \mu(\dd\varphi)
\]
is proportional to $\cos(n\psi)$
if and only if $\mu$ has a uniform density on the interval $[0,\pi]$.
Extending $\mu$ to the interval $[0,2\pi]$ as a symmetric (with respect
to the center $\pi$) measure, $\mu(A)=\mu(2\pi-A)$,
and defining the measure $\tilde\mu$ as $\tilde\mu(A)=\mu(2A)/2$ for
all Borel sets $A\in[0,\pi]$, we obtain
\begin{eqnarray*}
&&\int_0^\pi{\ln(\cos\varphi-\cos
\psi)^2}\cos(n\varphi) \mu(\dd\varphi) \\
&&\qquad=\frac{1}{2}\int
_0^{2\pi} \cos(n\varphi)\ln(\cos\varphi-\cos
\psi)^2 \mu(\dd\varphi)
\\
&&\qquad=\frac{1}{2}\int_0^{2\pi} \cos(n
\varphi)\ln\biggl(2\sin\frac
{\varphi
-\psi}{2}\sin\frac{\varphi+\psi}{2}
\biggr)^2 \mu(\dd\varphi)
\\
&&\qquad=\frac{1}{2}\int_0^{2\pi} \cos(n
\varphi)\ln2^2\mu(\dd\varphi) \\
&&\qquad\quad{}+\frac{1}{2}\int
_0^{2\pi} \cos(n\varphi)\ln\biggl(\sin
\frac
{\varphi-\psi
}{2} \biggr)^2\mu(\dd\varphi)
\\
&&\qquad\quad{}+\frac{1}{2}\int_0^{2\pi}
\cos(n\varphi)\ln\biggl(\sin\frac
{\varphi
+\psi}{2} \biggr)^2\mu(\dd
\varphi)
\\
&&\qquad=0 +\int_0^{\pi} \cos(2n\varphi)\ln
\sin^2(\varphi-\psi/2)\tilde\mu(\dd\varphi)\\
&&\qquad\quad{} +\int
_0^{\pi} \cos(2n\varphi)\ln\sin^2(
\varphi+\psi/2)\tilde\mu(\dd\varphi)
\\
&&\qquad=2\int_0^{\pi} \cos(2n\varphi-n\psi+n
\psi)\ln\sin^2(\varphi-\psi/2)\tilde\mu(\dd\varphi)
\\
&&\qquad=2\cos(n\psi)\int_0^{\pi} \cos(2n
\varphi-n\psi)\ln\sin^2(\varphi-\psi/2)\tilde\mu(\dd\varphi)
\\
&&\qquad\quad{}+2\sin(n\psi)\int_0^{\pi} \sin(2n
\varphi-n\psi)\ln\sin^2(\varphi-\psi/2)\tilde\mu(\dd\varphi).
\end{eqnarray*}
The ``if'' part follows from the facts that the
functions $\cos(2nz)\ln\sin^2(z)$ and $\sin(2nz)\ln\sin^2(z)$ are
$\pi$-periodic
and
\begin{eqnarray*}
\int_0^{\pi} \sin(2n\varphi-n\psi)\ln
\sin^2(\varphi-\psi/2)\frac{\dd
\varphi}{\pi}&=& \int_0^{\pi}
\sin(2n\varphi)\ln\sin^2(\varphi)\frac{\dd\varphi
}{\pi}=0,
\\
\int_0^{\pi} \cos(2n\varphi-n\psi)\ln
\sin^2(\varphi-\psi/2)\frac{\dd
\varphi}{\pi}&=& \int_0^{\pi}
\cos(2n\varphi)\ln\sin^2(\varphi)\frac{\dd\varphi
}{\pi}=-1/n.
\end{eqnarray*}

To prove\vspace*{1pt} the ``only if'' part
we need to show that the convolution of $\cos(2nz)\*\ln\sin^2(z)$ and
$\tilde\mu(z)$, that is,
\[
\int_0^{\pi} \cos\bigl(2n(\varphi-t)\bigr)\ln
\sin^2(\varphi-t)\tilde\mu(\dd\varphi)
\]
is constant for almost all $t\in[0,\pi]$ if and only if $\tilde\mu$
is uniform;
and the same holds for the convolution of $ \sin(2nz)\ln\sin^2(z)$ and
$\tilde\mu(z)$.
This, however, \mbox{follows} from Schmidt and Zhigljavsky
[(\citeyear{schzhi2009}), Lemma 3]
since\break $\cos(2nz)\ln\sin^2(z)\in L^2([0,\pi])$, and all complex Fourier
coefficients of these functions are nonzero.
Indeed,
\begin{eqnarray*}
\int_0^{\pi}\cos(2nt)\ln\sin^2(t)
\sin(2kt)\,\dd t&=&0 \qquad\forall k\in\mathbb{Z},
\\
\int_0^{\pi}\cos(2nt)\ln\sin^2(t)
\cos(2kt)\,\dd t&=&(\gamma_{|n+k|}+\gamma_{|n-k|})/2 \qquad\forall k
\in\mathbb{Z},
\end{eqnarray*}
where $\gamma_0=-2\pi\log2$ and $\gamma_k=-\pi/k$ for $k\in
\mathbb{N}$;
see formula 4.384.3 in \citet{graryz1965}.
\end{pf*}
\end{appendix}

\section*{Acknowledgments}

Parts of this paper were written during a visit of the authors at the
Isaac Newton Institute, Cambridge, UK, and the authors would like to
thank the institute for its hospitality and financial support.
We are also grateful to the referees and the Associate Editor for their
constructive comments on earlier versions of this manuscript.



\printaddresses

\end{document}